\documentclass{article}

\usepackage{amssymb}
\usepackage{amsmath}
\usepackage[all]{xy}
\usepackage[frenchb]{babel}
\usepackage{a4wide}
\usepackage{enumerate,fancyhdr}

\begin{document}

\title{Points de torsion sur les vari\'et\'es ab\'eliennes de type $\GSp$}

\author{Marc Hindry\footnote{hindry@math.jussieu.fr}, Nicolas Ratazzi \footnote{nicolas.ratazzi@math.u-psud.fr}
}

\renewcommand{\epsilon}{\varepsilon}
\renewcommand{\tilde}{\widetilde}

\newcommand{\Lie}{\textnormal{Lie}}
\newcommand{\Gal}{\textnormal{Gal}}
\newcommand{\tors}{\textnormal{tors}}
\newcommand{\Ker}{\textnormal{Ker}}
\newcommand{\red}{\textnormal{r\'{e}d}}
\newcommand{\Z}{\mathbb{Z}}
\newcommand{\Q}{\mathbb{Q}}
\newcommand{\C}{\mathbb{C}}
\newcommand{\R}{\mathbb{R}}
\newcommand{\N}{\mathbb{N}}
\renewcommand{\P}{\mathbb{P}}
\newcommand{\F}{\mathbb{F}}
\newcommand{\G}{\mathbb{G}}
\newcommand{\SD}{\mathbb{S}}
\renewcommand{\hat}{\widehat}

\newcommand{\T}{\textnormal{T}}

\newcommand{\PSp}{\textnormal{PSp}}
\newcommand{\K}{\overline{K}}

\newcommand{\mP}{\mathcal{P}}
\newcommand{\mA}{\mathcal{A}}
\renewcommand{\O}{\mathcal{O}}
\renewcommand{\L}{\mathcal{L}}
\newcommand{\V}{\mathcal{V}}
\newcommand{\kk}{\mathbb{K}}

\newcommand{\im}{\textnormal{Im}}

\newcommand{\m}{\mathfrak{m}}
\newcommand{\I}{\mathfrak{I}}

\newcommand{\hl}{\mathfrak{h}_{\ell}}
\newcommand{\el}{\mathfrak{e}_{\ell}}
\newcommand{\gsp}{\mathfrak{gsp}}
\newcommand{\n}{\mathfrak{n}}
\newcommand{\Tr}{\textnormal{Tr}}

\newcommand{\VV}{\mathsf{V}}

\newcommand{\mmid}{\|\|}

\newcommand{\Proj}{\textnormal{Proj\,}}
\newcommand{\Spec}{\textnormal{Spec\,}}
\newcommand{\Deg}{\textnormal{Deg\,}}
\newcommand{\e}{\varepsilon}
\newcommand{\ab}{\textnormal{ab}}
\newcommand{\End}{\textnormal{End}}
\renewcommand{\H}{\textnormal{H}}
\newcommand{\MT}{\textnormal{MT}}
\newcommand{\codim}{\textnormal{codim\,}}
\newcommand{\Hom}{\textnormal{Hom}}
\newcommand{\Aut}{\textnormal{Aut}}
\newcommand{\Hdg}{\textnormal{Hdg}}
\newcommand{\GL}{\textnormal{GL}}
\newcommand{\SL}{\textnormal{SL}}
\newcommand{\SLG}{\textnormal{SLG}}
\newcommand{\Sp}{\textnormal{Sp}}
\newcommand{\ssp}{\mathfrak{sp}}
\newcommand{\PGSp}{\textnormal{PGSp}}

\newcommand{\GSp}{\textnormal{GSp}}
\newcommand{\Res}{\textnormal{Res}}
\newcommand{\tr}{\textnormal{Tr}}
\newcommand{\frob}{\textnormal{Frob}}
\renewcommand{\ss}{\sigma}
\newcommand{\pr}{\textnormal{pr}}
\newcommand{\hdg}{\textnormal{Hdg}}

\renewcommand{\text}{\textnormal}

\newcommand{\ee}{\overline{e}}
\newcommand{\ff}{\overline{f}}
\newcommand{\tf}{\hat{f}}
\newcommand{\te}{\hat{e}}

\newtheorem{theo}{Th{\'e}or{\`e}me} [section]
\newtheorem{lemme}[theo]{Lemme}
\newtheorem{conj}{Conjecture}[section]
\newtheorem{prop}[theo]{Proposition}
\newtheorem{cor}[theo]{Corollaire}

\newcommand{\defi}{\addtocounter{theo}{1}{\noindent \textbf{D{\'e}finition \thetheo\ }}}
\newcommand{\rem}{\addtocounter{theo}{1}{\noindent \textbf{Remarque \thetheo\ }}}

\newcommand{\demo}{\noindent \textit{D{\'e}monstration} : }
\newcommand{\findemo}{\hfill \Box}

\maketitle

\hrulefill

\bigskip

\noindent \textbf{R\'esum\'e : } 
Soit $A$ une vari\'et\'e ab\'elienne d\'efinie sur un corps de nombres $K$, le nombre de points de torsion
d\'efinis sur une extension finie $L$   est born\'e polynomialement en terme du degr\'e $[L:K]$. Lorsque $A$ est isog\`ene \`a un produit de vari\'et\'es ab\'eliennes  simples de type $\GSp$, c'est-\`a-dire dont le groupe de Mumford-Tate est ``g\'en\'erique" (isomorphe au groupe des similitudes symplectiques) et  v\'erifiant la conjecture de Mumford-Tate, nous calculons l'exposant optimal dans cette borne, en terme de la dimension des sous-vari\'et\'es ab\'eliennes de $A$. Le r\'esultat est inconditionnel pour un produit de vari\'et\'es ab\'eliennes   simples dont l'anneau d'endomorphismes est $\Z$ et la dimension n'appartient pas \`a un ensemble exceptionnel explicite $\mathcal{S}=\{4,10,16,32,\dots\}$. Par ailleurs nous prouvons, suivant une strat\'egie de Serre, que si la conjecture de Mumford-Tate   est vraie pour des vari\'et\'es ab\'eliennes   de type $\GSp$, alors la conjecture de Mumford-Tate  est vraie pour un produit de telles vari\'et\'es ab\'eliennes.
  
\bigskip

\hrulefill

\bigskip

\noindent \textbf{Abstract : }
 Let $A$ be an abelian variety defined over a number field $K$, the number of torsion points rational over a finite extension $L$ is bounded polynomially in terms of the degree $[L:K]$.  When $A$ is isogenous to a product of simple abelian varieties of $\GSp$ type, i.e. whose Mumford-Tate group is ``generic" (isomorphic to the group of symplectic similitudes) and which satisfy the Mumford-Tate conjecture, we compute the optimal exponent for this bound in terms of the dimensions of the abelian subvarieties of $A$. The result is unconditional for a product of simple abelian varieties with endomorphism ring $\Z$ and dimension outside an explicit exceptional set $\mathcal{S}=\{4,10,16,32,\dots\}$. Furthermore, following a strategy of Serre, we also prove that if the Mumford-Tate conjecture is true for some abelian varieties of $\GSp$ type, it is then true for a product of such abelian varieties.

\hrulefill

\section{Introduction}

Soit $A$ une vari\'et\'e ab\'elienne de dimension $g\geq 1$ d\'efinie sur un corps de nombres $K$ plong\'e dans $\C$. Apr\`es avoir choisi une polarisation, on sait que le groupe de Mumford-Tate de $A$ (dont la d\'efinition est rappel\'ee au paragraphe \ref{MTH}), not\'e $\MT(A)$, est un sous-groupe alg\'ebrique sur $\Q$ du groupe des similitudes symplectiques not\'e $\GSp_{2g}$. Nous dirons que ``$A$ est de type $\GSp$" si son groupe de Mumford-Tate est g\'en\'erique, c'est-\`a-dire si $\MT(A)=\GSp_{2g}$. Une condition n\'ecessaire est d'avoir $\End_{\overline{K}}A=\Z$; cette condition n'est pas en g\'en\'eral suffisante, mais on sait (voir par exemple \cite{pink1}) qu'elle l'est si $g$ n'appartient pas \`a l'ensemble exceptionnel $\mathcal{S}$ d\'efini comme suit.

\medskip

\noindent{\bf Notation.}  On note $\mathcal{S}$ l'ensemble des entiers $g\geq 1$ tels que $2g$ est une puissance $k$-i\`eme avec $k\geq 3$ impair ou soit de la forme ${2k \choose k}$ avec   $k\geq 3$ impair; en symboles~:
\begin{equation}\label{defdeS}
\mathcal{S}:=\left\{g\geq 1\;|\; \exists k\geq 3,\;{\rm impair},\; \exists a\geq 1,\; g=2^{k-1}a^k \right\}\cup
\left\{g\geq 1\;|\; \exists k\geq 3,\;{\rm impair},\; g=\frac{1}{2}{2k\choose k} \right\}
\end{equation}

\medskip
\noindent Introduisons maintenant l'invariant que nous allons \'etudier.
\medskip

\noindent{\bf D\'efinition.} On pose 
\[\gamma(A)=\inf\left\lbrace x>0\, | \,  \forall F/K_0 \text{ finie, }\ \left|A(F)_{\tors}\right|\ll [F:K]^x\right\rbrace.\]

\noindent  La notation $\ll$ signifie qu'il existe  une constante  $C$, ne d\'ependant que de $A/K$, telle que l'on a $\left|A(F)_{\tors}\right|\leq C [F:K]^x$. 
On peut traduire la d\'efinition en le fait que $\gamma(A)$ est l'exposant le plus petit possible tel que pour tout $\epsilon>0$, il existe une constante $C(\epsilon)=C(A/K,\epsilon)$ telle que pour toute extension finie $F/K$ on a
\[\left|A(F)_{\tors}\right|\leq C(\epsilon) [F:K]^{\gamma(A)+\epsilon}.\]
On voit facilement que l'invariant d\'efini ci-dessus est ind\'ependant du corps de d\'efinition $K$ choisi et ne d\'epend en fait que de la classe d'isog\'enie de la vari\'et\'e ab\'elienne $A$. Comme toute vari\'et\'e ab\'elienne est isog\`ene \`a un  produit de vari\'et\'es ab\'eliennes simples et  toute vari\'et\'e ab\'elienne, apr\`es \'eventuellement une extension finie des scalaires, est isog\`ene \`a une vari\'et\'e ab\'elienne principalement polaris\'ee, on pourra donc sans dommage  imposer les conditions suivantes.

\medskip

\noindent{\bf Convention.} On supposera que la vari\'et\'e ab\'elienne $A$ est isomorphe \`a un produit de vari\'et\'es ab\'eliennes principalement polaris\'ees $A_1^{n_1}\times\dots\times A_d^{n_d}$ (avec des $A_i$ simples non isog\`enes deux \`a deux) et qu'on a remplac\'e $K$ par une extension finie convenable.

 \medskip 
\noindent Un r\'esultat g\'en\'eral d\^u \`a Masser \cite{lettremasser} donne une borne simple~:
\[\gamma(A)\leq\dim A\ \ \]
Cette borne est optimale lorsque $A$ est une puissance d'une courbe elliptique avec multiplication complexe; il est fort probable que la borne de Masser n'est jamais optimale dans les autres cas. L'invariant $\gamma(A)$ est calcul\'e dans \cite{hindry-ratazzi1} pour un produit de courbes elliptiques et, de mani\`ere diff\'erente, dans \cite{ratazziens} pour une vari\'et\'e ab\'elienne de type CM. Le probl\`eme analogue pour les  modules de Drinfeld est trait\'e dans \cite{florian}. Ces calculs nous ont  amen\'e \`a poser la question suivante.

\medskip

\begin{conj}\label{conjR} Soient $d$ et $n_1,\ldots, n_d$ des entiers strictement positifs. Si $A$ est isog\`ene \`a un pro\-duit de vari\'et\'es ab\'eliennes $\prod_{i=1}^d A_i^{n_i}$ avec les  $A_i$ simples non isog\`enes deux \`a deux, alors 
\begin{equation}\label{valeurconjdegamma}
\gamma(A)=\max_{\emptyset\not= I\subset\{1,\ldots,n\}}\frac{2\sum_{i\in I}n_i\dim A_i}{\dim\MT\left(\prod_{i\in I}A_i\right)}.
\end{equation}Dans le cas particulier o\`u $A_i$ est dimension $g_i$ et de type $\GSp$ on devrait avoir~:
\[\gamma(A)=\max_{\emptyset\not= I\subset\{1,\ldots,n\}}\frac{2\sum_{i\in I}n_ig_i}{1+\sum_{i\in I}2g_i^2+g_i}.\]
\end{conj}

\noindent Il est assez \'el\'ementaire de voir que $\gamma(A)$ est toujours sup\'erieur ou \'egal au membre de droite de (\ref{valeurconjdegamma}), la preuve est donn\'ee dans \cite{hindry-ratazzi1} (proposition 1.5). Il est aussi facile de voir que cette conjecture repose sur le cas particulier de la conjecture de Mumford-Tate
suivant. Notons $\rho_{\ell}:\Gal(\overline{K}/K)\rightarrow\GL(\T_{\ell}(A))\simeq\GL_{2g}(\Z_{\ell})$ la repr\'esentation $\ell$-adique associ\'ee \`a $A$ (avec $\T_{\ell}(A)=\varprojlim A[\ell^n]$ le module de Tate $\ell$-adique de $A$) et $G_{\ell}$ son image. On supposera pour simplifier (cf. la convention ci-dessus) que $A$ est munie d'une polarisation principale.

\begin{conj}\label{conjMT}\textnormal{\textbf{(Mumford-Tate)}} Si $A/K$ est une vari\'et\'e ab\'elienne de dimension $g$ sur un corps de nombres, de type $\GSp$, alors 
$G_{\ell}$ est un sous-groupe de $\GSp_{2g}(\Z_{\ell})$ d'indice fini, born\'e ind\'e\-pen\-damment de $\ell$.
\end{conj}

\noindent{\bf Remarque.} L'\'enonc\'e donn\'e ici est une version l\'eg\`erement plus forte que la conjecture initiale formul\'ee
dans \cite{mumfav} qui stipulait seulement une \'egalit\'e d'alg\`ebres de Lie $\ell$-adiques, \'equivalente \`a  la commensurabilit\'e de $G_{\ell}$ et $\MT(\Z_{\ell})$. Cependant dans bon nombre de cas, on peut 
d\'emontrer que la conjecture initiale suffit \`a entra\^{\i}ner la conjecture plus forte o\`u l'indice fini 
d\'ependant {\it a priori} de $\ell$ est en fait uniform\'ement born\'e; c'est notamment le cas pour les 
vari\'et\'es ab\'eliennes de type $\GSp$ (voir ci-dessous).

\noindent Dans la direction de la conjecture de Mumford-Tate pour les vari\'et\'es ab\'eliennes, on dispose de divers r\'esultats \cite{bgk1,bgk2,chi,hall,pink1,ribetrm,college8485}; dans notre contexte le r\'esultat important est  un  th\'eor\`eme de Serre \cite{college8485} compl\'et\'e par Pink \cite{pink1}, o\`u  l'ensemble $\mathcal{S}$ est d\'ecrit par l'\'equation (\ref{defdeS}) ci-dessus, et dans une autre direction par Hall \cite{hall}.

\begin{theo}\label{th1}\textnormal{\textbf{(Serre, Pink, Hall)}} Si $A/K$ est une vari\'et\'e ab\'elienne de dimension $g$ n'appartenant pas \`a $\mathcal{S}$, d\'efinie sur un corps de nombres, de type $\GSp$ alors 
$G_{\ell}$ est un sous-groupe de $\GSp_{2g}(\Z_{\ell})$ d'indice fini, born\'e ind\'e\-pen\-damment de $\ell$. Si $g$ est quelconque mais l'on suppose que le mod\`ele de N\'eron de $A$ sur $\mathcal{O}_K$ poss\`ede une fibre semistable avec dimension torique \'egale \`a un,  la m\^eme conclusion vaut.
\end{theo}
\medskip Dans \cite{serremfv} (th\'eor\`eme 3 paragraphe 7.), Serre d\'emontre  que, dans le cas de type $\GSp$, le fait que le groupe $G_{\ell}$ est un sous-groupe de $\GSp_{2g}(\Z_{\ell})$ d'indice fini entra\^{\i}ne que $G_{\ell}$ est {\em  \'egal} \`a 
 $\GSp_{2g}(\Z_{\ell})$ pour $\ell$ assez grand et que cela est vrai si $g$ est impair ou valant $2$ ou $6$ ; dans \cite{pink1}, Pink  \'etablit  que le groupe $G_{\ell}$ est un sous-groupe de $\GSp_{2g}(\Z_{\ell})$ d'indice fini pour les valeurs de $g$  \'evitant l'ensemble $\mathcal{S}$. On remarquera que, si $g\notin\mathcal{S}$, la vari\'et\'e ab\'elienne $A$ est de type GSp si et seulement si $\End_{\bar{K}}(A)=\Z$. Enfin Hall \cite{hall} montre que, si l'on suppose que le mod\`ele de N\'eron de $A$ sur $\mathcal{O}_K$ poss\`ede une fibre semistable avec dimension torique \'egale \`a $1$, alors on obtient la m\^eme conclusion.

\begin{theo}\label{th2}  Si $A/K$ est une vari\'et\'e ab\'elienne de dimension $g$ telle que, pour tout premier $\ell$,  le groupe de Galois associ\'e $G_{\ell}$ est d'indice fini  dans  $\GSp_{2g}(\Z_{\ell})$, alors  
\begin{equation}\label{valeurgamma}
\gamma(A)=\frac{2\dim A}{\dim\MT(A)}=\frac{2g}{2g^2+g+1}.
\end{equation}
\end{theo}

\begin{cor}\label{coroth2}  Si $A/K$ est une vari\'et\'e ab\'elienne d\'efinie sur un corps de nombres $K$, de dimension $g$ telle que $\End_{\overline{K}}A=\Z$. Supposons l'une des deux conditions suivantes r\'ealis\'ee:
\begin{enumerate}
\item la dimension $g$ n'appartient pas \`a $\mathcal{S}$,
\item le mod\`ele de N\'eron de $A$ sur $\mathcal{O}_K$ poss\`ede une fibre semistable avec dimension torique \'egale \`a $1$,
\end{enumerate} 
alors
 \begin{equation}\label{valeurgammacor}
\gamma(A)=\frac{2\dim A}{\dim\MT(A)}=\frac{2g}{2g^2+g+1}.
\end{equation}
\end{cor}

\medskip

\noindent Nous voulons ensuite \'etendre ce r\'esultat au produit de vari\'et\'es ab\'eliennes du m\^eme type. Pour cela nous d\'emontrons \'egalement que la conjecture de Mumford-Tate forte est vraie pour un tel produit de vari\'et\'es ab\'eliennes. Nous renvoyons au paragraphe \ref{MTH} pour les d\'efinitions des groupes de Hodge et de Mumford-Tate.

\begin{theo}\label{prophodgeprod} Soient $r$ et $n_1\ldots,n_r$ des entiers strictement positifs. Soient $A_i$ des vari\'et\'es ab\'e\-liennes de dimension $g_i$ non isog\`enes deux \`a deux telles que $\hdg(A_i)=\Sp_{2g_i}$. Posons $A:=A_1^{n_1}\times\dots\times A_r^{n_r}$ et, pour tout premier $\ell$, notons $\rho_{\ell,i}$ (respectivement $\rho_{\ell}=\rho_{\ell,1}\times\ldots,\rho_{\ell,r}$) les repr\'esentations $\ell$-adiques associ\'ees aux $A_i$ (respectivement \`a $A$) alors~:
\begin{enumerate}
\item l'inclusion naturelle suivante est un isomorphisme~:
$$\hdg(A)\cong\hdg\left(A_1\times\dots\times A_r\right)\hookrightarrow\Sp_{2g_1}\times\dots\times\Sp_{2g_r}.$$
\item soit $\ell$ un nombre premier. Si  pour tout $i$, on a $\rho_{\ell,i}\left(\Gal(K(A_i[\ell^{\infty}])/K(\mu_{\ell^{\infty}}))\right)\cong \Sp_{2g_i}(\Z_{\ell})$ (\`a indice fini pr\`es) 
alors, on a (\`a indice fini pr\`es)~:
$$ \rho_{\ell}\left(\Gal(K(A[\ell^{\infty}])/K(\mu_{\ell^{\infty}}))\right)\cong \Sp_{2g_1}(\Z_{\ell})\times\dots\times\Sp_{2g_r}(\Z_{\ell}).$$ 
\item si de plus l'indice fini pour chaque $A_i$ est born\'e ind\'ependamment de $\ell$, il en est de m\^eme pour $A$.
\end{enumerate}
\end{theo}
\rem On peut r\'esumer le th\'eor\`eme en disant que si la conjecture de Mumford-Tate   est vraie pour $A_i$ de type $\GSp$, alors la conjecture de Mumford-Tate  est vraie pour $A=\prod_iA_i^{n_i}$. Sous la derni\`ere hypoth\`ese (point 3.), on peut montrer que (et nous allons le faire dans la preuve), il y a en fait \'egalit\'e pour tout $\ell$ suffisamment grand.

\medskip

\noindent Nous pouvons ainsi obtenir la valeur de $\gamma(A)$ pour un produit de vari\'et\'es ab\'eliennes de type $\GSp$ v\'erifiant Mumford-Tate fort.

\begin{theo}\label{thproduit}  Soit $A_i/K$ des vari\'et\'es ab\'eliennes non isog\`enes deux \`a deux, de dimension $g_i$ d\'efinies sur un corps de nombres, de type $\GSp$ telles que, pour tout premier $\ell$, le groupe de Galois associ\'e $G_{\ell}$ est d'indice fini   dans  $\GSp_{2g_i}(\Z_{\ell})$. Soit 
$A:=A_1^{n_1}\times\dots\times A_r^{n_r}$ avec des entiers   $n_i\geq 1$. On a alors
\begin{equation}\label{valeurgammaproduit}
\gamma(A)=\max_{\emptyset\not= I\subset\{1,\ldots,n\}}\left\{\frac{2\sum_{i\in I}n_i\dim A_i}{\dim\MT(\prod_{i\in I}A_i)}\right\}=\max_{\emptyset\not= I\subset\{1,\ldots,n\}}\left\{\frac{2\sum_{i\in I}n_ig_i}{1+\sum_{i\in I}2g_i^2+g_i}\right\}.
\end{equation}
\end{theo}

\begin{cor}\label{corothproduit}  Soit $A_i/K$ des vari\'et\'es ab\'eliennes non isog\`enes deux \`a deux, de dimension $g_i$ d\'efinies sur un corps de nombres, telles que $\End_{\overline{K}}A_i=\Z$. Soit 
$A:=A_1^{n_1}\times\dots\times A_r^{n_r}$ avec des entiers   $n_i\geq 1$. Supposons que pour chaque $A_i$ l'une des deux conditions suivantes soit r\'ealis\'ee:
\begin{enumerate}
\item la dimension $g_i$ n'appartient pas \`a $\mathcal{S}$,
\item le mod\`ele de N\'eron de $A_i$ sur $\mathcal{O}_K$ poss\`ede une fibre semistable avec dimension torique \'egale \`a un,
\end{enumerate} 
alors
\begin{equation}\label{valeurgammaproduitcor}
\gamma(A)=\max_{\emptyset\not= I\subset\{1,\ldots,n\}}\left\{\frac{2\sum_{i\in I}n_i\dim A_i}{\dim\MT(\prod_{i\in I}A_i)}\right\}=\max_{\emptyset\not= I\subset\{1,\ldots,n\}}\left\{\frac{2\sum_{i\in I}n_ig_i}{1+\sum_{i\in I}2g_i^2+g_i}\right\}.
\end{equation}
\end{cor}

\medskip

\noindent L'organisation du texte est la suivante. Le deuxi\`eme paragraphe contient divers rappels et r\'esultats de th\'eorie des groupes, pour la plupart \'el\'ementaires concernant les groupes symplectiques et le cardinal de l'ensemble des points d'un sous-groupe alg\'ebrique de $\GL_n$ \`a valeurs dans $\Z/\ell^m\Z$. Le troisi\`eme paragraphe d\'ecrit quelques sous-modules du module de Tate $\ell$-adique d'une vari\'et\'e ab\'elienne et aborde le calcul de l'intersection de l'extension de $K$ engendr\'ee par des points d'ordre fini avec l'extension cyclotomique. Le quatri\`eme paragraphe contient la preuve du th\'eor\`eme \ref{th2} (voir le th\'eor\`eme \ref{th3}), i.e. le calcul de l'invariant $\gamma(A)$ pour une vari\'et\'e ab\'elienne simple g\'en\'erique.  Le cinqui\`eme paragraphe contient des rappels sur les groupes de Mumford-Tate et de Hodge, et on y calcule ces groupes pour un produit de vari\'et\'es ab\'eliennes simples g\'en\'eriques.  Le sixi\`eme paragraph!
 e contient la preuve du th\'eor\`eme \ref{thproduit}, i.e. le calcul de l'invariant $\gamma(A)$ pour un produit de vari\'et\'es ab\'eliennes simples g\'en\'eriques.

\medskip

\noindent{\bf Remerciements.}  Les deux auteurs remercient le referee pour ses remarques et corrections pertinentes.

\section{Rappels et lemmes de groupes} 

\medskip
Nous rassemblons dans ce paragraphe des lemmes combinatoires, des ``lemmes de comptage" de points \`a valeurs dans $\Z/\ell^n\Z $ de divers  groupes alg\'ebriques, des lemmes sp\'ecifiques aux groupes $\Sp_{2n}$ et enfin une description de certains stabilisateurs sous-groupes de $\Sp_{2n}$ qui seront importants pour la d\'emonstration du th\'eor\`eme \ref{th2}.

\medskip 

\noindent  Dans toute la suite, la notation $\gg\ll$ signifiera ``\`a une constante multiplicative pr\`es, ind\'ependante de $\ell$''.

\subsection{Lemmes de comptages}

\begin{lemme}\label{hensel}Soit $G/\Z_{\ell}$ un sous-groupe alg\'ebrique de $\GL_n$, de dimension $d$, tel que la r\'eduction de $G$ sur $\mathbb{F}_{\ell}$ est un groupe lisse. On a 
\[\forall{m\geq 1},\ \ \ \left|G\left(\Z/\ell^m\Z\right)\right|=\ell^{(m-1)d}\left|G\left(\Z/\ell\Z\right)\right|.\]
\end{lemme}
\demo Il s'agit d'une variante du lemme de Hensel. Plus pr\'ecis\'ement, consid\'erons l'application naturelle (pour $m\geq 1$):
$$\pi_m: G(\Z/\ell^{m+1}\Z)\rightarrow G(\Z/\ell^m\Z).$$
L'hypoth\`ese de lissit\'e entra\^{\i}ne d'une part que  $\pi_m$ est surjective et d'autre part que le noyau s'identifie avec l'espace tangent \`a $G$ sur $\F_{\ell}$ et a donc pour ordre $\ell^{\dim G}=\ell^d$.\hfill $\Box$

\medskip

\begin{prop}\label{comptage} Soit $G/\Z$ un sous-groupe alg\'ebrique de $\GL_n$. Notons $d$ sa dimension (sur $\Q$), $r$ son rang et $n_G$ son nombre de composantes connexes. On a l'encadrement~:
\[C_{1,\ell}(1-\frac{1}{\ell})^r\leq \frac{|G(\Z/\ell\Z)|}{\ell^d}\leq n_GC_{2,\ell}(1+\frac{1}{\ell})^r\]
\noindent o\`u $C_1:=\prod_{\ell}C_{1,\ell}$ et $C_2:=\prod_{\ell}C_{2,\ell}$ sont des produits convergents.
Plus g\'en\'eralement on a~:
\[\forall N\geq 1,\ \ \ C_{1}^{\omega(N)}\prod_{\ell\,|\,N}(1-\frac{1}{\ell})^r\leq \frac{|G(\Z/N\Z)|}{N^d}\leq C_{2}^{\omega(N)}\prod_{\ell\,|\,N}(1+\frac{1}{\ell})^r,\]
o\`u $\omega (x)$ est le nombre de premiers divisants $x$.
\end{prop}
\demo Il s'agit d'un r\'esultat qui est sans doute bien connu des experts, mais dont nous pr\'ef\'erons redonner une preuve rapide ici. Notons tout d'abord que la seconde partie de la proposition d\'ecoule directement de la premi\`ere et du lemme \ref{hensel} pr\'ec\'edent. Il s'agit donc de mesurer le cardinal de $G(\Z/\ell\Z)$. Dans toute la discussion qui suit il convient de traiter \`a part un petit nombre (fini !) de premiers. Ceci ne pose pas de probl\`eme : pour tout $\ell$, il existe $C_{1,\ell}, C_{2,\ell}>0$  telles que 
\[C_{1,\ell}\leq \frac{|G(\Z/\ell\Z)|}{\ell^d}\leq C_{2,\ell}.\]
\noindent Dans la suite nous nous autoriserons donc toujours \`a supposer que $\ell$ est pris en dehors d'un ensemble fini (d\'ependant du groupe $G$). Par ailleurs, si le groupe $G$ n'est pas connexe, notons $n_G$ le nombre de composantes connexes et $G^0$ la composante connexe de l'identit\'e de $G$. On a l'encadrement 
\[|G^0(\Z/\ell\Z)|\leq |G(\Z/\ell\Z)|\leq n_G|G^0(\Z/\ell\Z)|.\]
\noindent Dans la suite, nous pouvons donc supposer que $G$ est connexe. Sur $\Q$ le groupe $G$ se d\'ecompose comme produit semi-direct sous la forme $G=UR$ avec $U$ unipotent et $R$ r\'eductif (il s'agit de la d\'ecomposition de L\'evi).    
On a donc
\[|G(\Z/\ell\Z)|=|U(\Z/\ell\Z)|\times |R(\Z/\ell\Z)|.\]
\noindent Or pour un groupe unipotent, l'exponentielle donne un isomorphisme entre $U$ et son alg\`ebre de Lie. Notamment, on a
\[|U(\Z/\ell\Z)|=\left|\mathbb{A}^{\dim U}(\Z/\ell\Z)\right|=\ell^{\dim U}.\]
\noindent Il suffit donc de prouver le r\'esultat pour un groupe r\'eductif $R$. Or ce dernier est $\Q$-isog\`ene au produit direct $T\times S$ o\`u $T$ est un tore et $S$ est semi-simple. En particulier ces groupes ont m\^eme nombres de points sur le corps $\Z/\ell\Z$ : 
\[|R(\Z/\ell\Z)|=|T(\Z/\ell\Z)|\times |S(\Z/\ell\Z)|.\]
\noindent Pour un groupe semi-simple, nous disposons d'un th\'eor\`eme de Chevalley (cf. \cite{steinberg}) assurant que le produit
\[\prod_{\ell}\frac{|S(\Z/\ell\Z)|}{\ell^{\dim S}}\]
\noindent converge. Il suffit donc de v\'erifier que le r\'esultat annonc\'e est vrai dans le cas d'un tore de dimension $r\geq 1$. Dans ce dernier cas, on dispose d'une formule exacte pour le cardinal de $T(\Z/\ell\Z)$ (cf. par exemple \cite{vosk} p.104 theorem 2) : 
\[|T(\Z/\ell\Z)|=\det(\ell I_r-h(\sigma)),\]
\noindent o\`u $\sigma$ est un g\'en\'erateur topologique de $\mathcal{G}_{\ell}:=\Gal(\overline{\mathbb{F}}_{\ell}/\mathbb{F}_{\ell})$ et o\`u $h : \mathcal{G}_{\ell}\rightarrow \GL_r$ est la repr\'esentation d\'efinie par le $\mathcal{G}_{\ell}$-module $X^*(T_{\mathbb{F}_{\ell}})$. Les valeurs propres de $h(\sigma)$ sont de module $1$, donc
\[(1-\frac{1}{\ell})^r\leq \frac{|T(\Z/\ell\Z)|}{\ell^d}\leq (1+\frac{1}{\ell})^r.\]
\noindent \hfill $\Box$

\medskip

\begin{cor}\label{compter} Soient $\mathcal{G}_1$ un sous-groupe alg\'ebrique sur $\Z$ de $\GL$ et soit $G_1$ un sous-groupe alg\'ebrique de $\GL_{\Z_{\ell}}$ sur $\Z_{\ell}$ tel que sa r\'eduction modulo $\ell$ est conjugu\'ee sur $\F_{\ell}$ \`a $\mathcal{G}_{1,\F_{\ell}}$. On a
\[
\forall m\geq 1,\ \  \forall\ell \ \  \text{tel que $\mathcal{G}_{1,\F_{\ell}}$ est lisse, on a  :}\ \left|G_1(\Z/\ell^m\Z)\right|\gg\ll \ell^{m\dim \mathcal{G}_1}.\]
\end{cor}
\demo Si $\ell$ est tel que $\mathcal{G}_{1,\F_{\ell}}$ est lisse, l'assertion est une cons\'equence directe du lemme \ref{hensel} et de la proposition \ref{comptage}. $\findemo$

\medskip

\begin{lemme}\label{cle} Soit $G$ un sous-groupe alg\'ebrique sur $\Z$ de $\GL$,  
soit $t\in\N^*$ et soit $\mathcal{G}_1,\ldots,\mathcal{G}_t$ une suite de sous-groupes alg\'ebriques de $G$ sur $\Z$.  Soient $G_1\subset G_{2}\subset\dots\subset G_t$ une suite de sous-groupes alg\'ebriques sur $\Z_{\ell}$ de $G_{\Z_{\ell}}$. On suppose que pour tout $i$, $G_i$ est conjugu\'e sur $\F_{\ell}$ \`a $\mathcal{G}_i$. On note $g_i:=\dim\mathcal{G}_i=\dim G_i$ et $d_i:=\codim_G \mathcal{G}_i=\codim_{G_{\Z_{\ell}}}G_i$ et on pose, pour toute suite croissante d'entiers   
$0=m_0<m_1<m_2<\dots<m_t$~:
\[H(m_1,\dots,m_t)=\left\{M\in G(\Z_{\ell})\;|\; M\in G_i\mod \ell^{m_i}\right\}.\]
\noindent Pour tous les $\ell$ tels que $G$ et les $\mathcal{G}_i$ sont lisses sur $\F_{\ell}$, on a alors
\[\left(G(\Z_{\ell}):H(m_1,\dots,m_t)\right)\gg\ll\ell^{\sum_{i=1}^td_i(m_i-m_{i-1})},\]
\noindent les constantes multiplicatives qui interviennent dans $\gg\ll$ ne d\'ependant que des $\mathcal{G}_i$ et pas des $G_i$.
\end{lemme}

\demo Pour $t=1$, consid\'erons l'homomorphisme de r\'eduction ${\rm r\acute{e}d}:G(\Z_{\ell})\rightarrow G(\Z/\ell^{m_1}\Z)$. Il est surjectif par l'hypoth\`ese de lissit\'e et il induit donc un isomorphisme entre $G(\Z_{\ell})/H(m_1)$ et $G(\Z/\ell^{m_1}\Z)/G_1(\Z/\ell^{m_1}\Z)$. En appliquant le corollaire \ref{compter} \`a $\mathcal{G}_1$ et $G_1$, on voit que : 
\[\left|G(\Z_{\ell})/H(m_1)\right|\gg\ll\ell^{m_1(\dim G-\dim \mathcal{G}_1)}=\ell^{m_1d_1}.\]
\noindent Montrons maintenant l'\'enonc\'e par induction sur $t$. Introduisons les groupes finis~:
$$G(m_1,\dots,m_t)=\{M\in G_t(\Z/\ell^{m_{t}}\Z)\;|\;M\in G_i\mod\ell^{m_i}\;({\rm pour}\; 1\leq i\leq t-1)\}.$$
L'\'enonc\'e qu'on veut d\'emontrer peut se traduire par 
$$(G(\Z/\ell^{m_t}\Z):G(m_1,\dots,m_t))\gg\ll\ell^{\sum_{i=1}^td_i(m_i-m_{i-1})}.$$
On regarde l'homomorphisme naturel $\phi: G_t(\Z/\ell^{m_t}\Z)\rightarrow G_t(\Z/\ell^{m_{t-1}}\Z)$; il est surjectif par l'hypoth\`ese de lissit\'e (sur $\F_{\ell}$) et son noyau est d'ordre $\ell^{g_t(m_t-m_{t-1})}$. En notant $H:=G(m_1,\dots,m_{t-1})$ le sous-groupe de $G_t(\Z/\ell^{m_{t-1}}\Z)$, on voit que $H'=\phi^{-1}(H)=G(m_1,\dots,m_t)$ est de cardinal
$\ell^{g_t(m_t-m_{t-1})}|G(m_1,\dots,m_{t-1})|$. On obtient ainsi
$$(G(\Z/\ell^{m_t}\Z):G(m_1,\dots,m_t))= \ell^{-g_t(m_t-m_{t-1})}\frac{|G(\Z/\ell^{m_t}\Z)|}{|G(\Z/\ell^{m_{t-1}}\Z)|}   (G(\Z/\ell^{m_t}\Z):G(m_1,\dots,m_{t-1})),$$
d'o\`u en appliquant l'hypoth\`ese de r\'ecurrence au dernier terme du membre de droite
\[(G(\Z/\ell^{m_t}\Z):G(m_1,\dots,m_t))\gg\ll
\ell^{d_t(m_t-m_{t-1})+\sum_{i=1}^{t-1}d_i(m_i-m_{i-1})}=\ell^{\sum_{i=1}^td_i(m_i-m_{i-1})}.\]
\noindent $\findemo$

\medskip

Nous incluons les deux \'enonc\'es \'el\'ementaires de th\'eorie des groupes suivants  pour future r\'ef\'erence~:

\begin{lemme}\label{lemgs} Soit $H$ un sous-groupe d'indice fini $m$ dans $G$, \'ecrivons
$G=\cup_{i=1}^mg_iH$ alors $\tilde{H}:=\cap_{i=1}^mg_iHg_i^{-1}$ est un sous-groupe normal d'indice divisant $m!$. En particulier, si $G$ est simple et $H$ est  d'un sous-groupe strict d'indice $m$ dans $G$, alors $m!\geq |G|$.
\end{lemme}
\demo  Consid\'erons l'action de $G$ sur $G/H$ donn\'ee par $(g,aH)\mapsto gaH$; le noyau de l'action est $\tilde{H}$, ce qui fournit un homomorphisme injectif $G/\tilde{H}\rightarrow\mathfrak{S}_m$.
$\findemo$

\medskip

\begin{lemme}\label{aut1} 
 Soit $H$ un sous-groupe normal de $G$ dont le centralisateur est trivial (i.e. tel que $C_G(H)=\{1\}$). Un automorphisme de $G$ induisant l'identit\'e sur $H$ est l'identit\'e sur $G$.
\end{lemme}  

\demo  Soit $\psi\in\Aut(G)$ tel que $\psi_{|H}=id_H$. Soit $x\in H$ et $y\in G$ alors
$\psi(y)x\psi(y)^{-1}=\psi(yxy^{-1})=yxy^{-1}$. On en tire $y^{-1}\psi(y)\in C_G(H)$ donc $\psi(y)=y$ et on a bien $\psi=id_G$.  
$\findemo$

\medskip

\noindent Enfin nous utiliserons le lemme combinatoire \'el\'ementaire suivant.

\begin{lemme}\label{combielem} Soit $a_1,\ldots,a_k$ et $b_1,\ldots,b_k$ des entiers strictement positifs. On a l'\'egalit\'e
\begin{equation}\label{formulelem}
\sup_{m_1\geq\dots\geq m_k}\left\{\frac{\sum_{i=1}^ka_im_i}{\sum_{i=1}^kb_im_i}\right\}=
\max_{1\leq h\leq k}\left\{\frac{\sum_{i=1}^ha_i}{\sum_{i=1}^hb_i}\right\}.
\end{equation}
\end{lemme}
\demo C'est le lemme 7.10 de \cite{ratazziens}. Il s'agit d'une simple application de la transformation d'Abel.$\findemo$

\medskip
\noindent Dans la version pour les produits de vari\'et\'es ab\'eliennes de type GSp nous utiliserons la g\'en\'erali\-sa\-tion (facile) suivante~:

\begin{lemme}\label{combielem2} Soit $d\geq 1$ un entier, et pour tout $i\in\{1,\ldots,d\}$, soit $t_i\geq 1$ des entiers. Pour $i\leq d$ et $j\leq t_i$, on se donne \'egalement des entiers $a_{ij}$ et $b_{ij}$, strictement positifs. On a l'\'egalit\'e
\begin{equation}\label{formulelem2}
\sup_{\underset{1\leq i\leq d}{m_{i1}\geq\dots\geq m_{i t_i}}}\left\{\frac{\sum_{i=1}^d\sum_{j=1}^{t_i}a_{ij}m_{ij}}{\sum_{i=1}^d\sum_{j=1}^{t_i}b_{ij}m_{ij}}\right\}=
\max_{\underset{1\leq i\leq d}{1\leq h_i\leq t_i}}\left\{\frac{\sum_{i=1}^d\sum_{j=1}^{h_i}a_{ij}}{\sum_{i=1}^d\sum_{j=1}^{h_i}b_{ij}}\right\}.
\end{equation}
\end{lemme}
\demo  Tout d'abord, en choisissant $1=m_{i1}=\dots=m_{ih}$ et $m_{ij}=0$ pour $j>h$, on voit que le membre de gauche est sup\'erieur ou \'egal au membre de doite. Pour obtenir l'in\'egalit\'e inverse, posons $t=\max t_i$ et posons $m_{ij}=a_{ij}=b_{ij}=0$ si $j\geq k_i+1$. Par transformation d'Abel, on a
\begin{align*} A        & := \sum_{i=1}^d\sum_{j=1}^{t_i}a_{ij}m_{ij}=\sum_{i=1}^d\sum_{j=1}^{t_i}\sum_{l=1}^ja_{il}(m_{ij}-m_{ij+1})\\
                                                                        & = \sum_{j=1}^t\sum_{i=1}^{d}\sum_{l=1}^ja_{il}(m_{ij}-m_{ij+1})\\
                                                                        & \leq \sum_{j=1}^t  \sup_{\underset{i\in \{1,\ldots, d\}}{1\leq j\leq t_i}}\left(\frac{\sum_{i=1}^d\sum_{l=1}^ja_{il}}{\sum_{i=1}^d\sum_{l=1}^jb_{il}}\right)                         \sum_{i=1}^{d}\sum_{l=1}^jb_{il}(m_{ij}-m_{ij+1})\\
                                                                        & \leq \sum_{j=1}^t  \sup_{\underset{i\in \{1,\ldots, d\}}{1\leq j\leq t_i}}\left(\frac{\sum_{i=1}^d\sum_{l=1}^ja_{il}}{\sum_{i=1}^d\sum_{l=1}^jb_{il}}\right)                         \sum_{i=1}^{d}b_{ij}m_{ij}\\
                                                                & \leq \sup_{\underset{i\in \{1,\ldots, d\}}{1\leq j\leq t_i}}\left(\frac{\sum_{i=1}^d\sum_{l=1}^ja_{il}}{\sum_{i=1}^d\sum_{l=1}^jb_{il}}\right)                   \sum_{i=1}^{d}\sum_{j=1}^{t_i}b_{ij}m_{ij}
\end{align*}
\noindent On conclut en divisant par la quantit\'e $ \sum_{i=1}^{d}\sum_{j=1}^{t_i}b_{ij}m_{ij}$. $\findemo$

\subsection{Les groupes $\Sp$ et $\GSp$\label{prs}}

\noindent Soit $J$ une matrice antisym\'etrique non d\'eg\'en\'er\'ee, on d\'efinit le groupe alg\'ebrique~:
$$\GSp_{2g,J}:=\left\{M\in\GL_{2g}\;|\; \exists \lambda(M)\in\G_m,\;^t\!{M}JM=\lambda(M) J\right\}.$$
Apr\`es changement de base, on peut supposer que $J=\begin{pmatrix} 0 & I_g\cr -I_g&0\cr\end{pmatrix}$; on note alors $\GSp_{2g}=\GSp_{2g,J}$. C'est un groupe alg\'ebrique  sur $\Z$,
 et on a 
\[M=\begin{pmatrix} A& B\cr C& D\cr\end{pmatrix}\in\GSp_{2g}\ \ \iff\ \  \left\{\begin{matrix} ^t\!AC\;\text{et}\; ^t\!BD\;\text{sont sym\'etriques}\cr  \exists \lambda(M)\in\G_m,\;^t\!AD-\,{}^t\!CB=\lambda(M) I_g\cr\end{matrix}\right.\]
\noindent On introduit $\lambda~:~\GSp_{2g}\rightarrow\G_m$, l'homomorphisme qui associe \`a $M$ son multiplicateur $\lambda(M)$.

\medskip

\rem\label{detlambda} Notons le lien suivant entre l'application multiplicateur et le d\'eterminant~:
\[\forall M\in\GSp_{2g}\ \ \  (\det M)=\lambda(M)^g.\]

\medskip

\noindent Par d\'efinition de $\lambda$, on a $\Ker\lambda=\Sp_{2g}$. Comme $\GSp_{2g}$ (ainsi que $\Sp_{2g}$) est stable par transposition on peut aussi en d\'eduire que $A\,^t\!B$ et $C\,^t\!D$ sont sym\'etriques.

\noindent L'alg\`ebre de Lie de $\Sp_{2g}$ s'identifie \`a l'alg\`ebre des matrices $M$ telles que
$^t\!MJ+JM=0$, c'est-\`a-dire \`a~:
$$\ssp_{2g}=\left\{M=\begin{pmatrix} A& B\cr C& D\cr\end{pmatrix}\ |\ D=-\,^t\!A\,\text{ et $B,C$ sont sym\'etriques}\right\}.$$
On v\'erifie ais\'ement sur l'alg\`ebre de Lie que $\dim\Sp_{2g}=\dim\ssp_{2g}=2g^2+g$.

\medskip

\noindent  Si $L_0$ d\'esigne le lagrangien (sous-espace isotrope de dimension maximale) engendr\'e par les $g$ premiers vecteurs de la base canonique, son fixateur dans $\Sp_{2g}$ est
$\left\{M=\begin{pmatrix}I_g&S\cr 0&I_g\cr\end{pmatrix}\ |\ S\text{ sym\'etrique }\right\}$ et son stabilisateur dans $\Sp_{2g}$ est 
$\left\{M=\begin{pmatrix}A&B\cr 0&\,^t\!A^{-1}\cr\end{pmatrix}\ |\ A\in\GL_g\text{ et }B\,^t\!A\text{ sym\'etrique }\right\}$.  Son fixateur dans $\GSp_{2g}$ est
$\left\{M=\begin{pmatrix}I_g&S\cr 0&\lambda I_g\cr\end{pmatrix}\ |\  S\text{ sym\'etrique et }\lambda\in\G_m\right\}$ et son stabilisateur dans le groupe $\GSp_{2g}$ est 
$\left\{M=\begin{pmatrix}A&B\cr 0&\,\lambda\,^t\!A^{-1}\cr\end{pmatrix}\ |\ A\in\GL_g\text{ et }B\,^t\!A\text{ sym\'etrique  et }\lambda\in\G_m\right\}$.

\medskip

\begin{lemme}\label{lemcong} Soit $M=\begin{pmatrix} A& B\\ C& D\end{pmatrix}\in\GSp_{2g}(\Z_{\ell})$ telle que
\[Me_1\equiv e_1\mod\ell^m\text{ et }
Me_{g+1}\equiv e_{g+1} \mod\ell^m,\]
alors $\lambda(M)\equiv 1\mod\ell^m$.
\end{lemme}
\demo Notons $\epsilon_1$ le vecteur colonne  \`a $g$ lignes et coordonn\'ees $1,0,\dots, 0$. Les hypoth\`eses se traduisent par
$$A\epsilon_1\equiv \epsilon_1\mod\ell^m,\; De_1\equiv \epsilon_1\mod\ell^m,\; C\epsilon_1\equiv 0\mod\ell^m\;\text{et}\; B\epsilon_1\equiv 0\mod\ell^m.$$
On en tire donc $\lambda(M) \epsilon_1=(\,^t\!AD-\,^t\!CB)(\epsilon_1)\equiv \begin{pmatrix}a_{11}\\ \vdots \\ a_{1g}\end{pmatrix}\mod\ell^m$ et donc $\lambda(M)\equiv a_{11}\mod\ell^m$. Comme par ailleurs
$\epsilon_1\equiv A\epsilon_1\equiv \begin{pmatrix}a_{11}\\ \vdots \\ a_{g1}\end{pmatrix}\mod\ell^m$, on a aussi $a_{11}\equiv 1\mod\ell^m$.
$\findemo$

\medskip

\begin{cor}\label{corB}  Pour tout entier $m$, soit $G(m)$ le sous-groupe
\[G(m):=\left\{M\in\GSp_{2g}(\Z_{\ell})\ |\ \forall x,\ Mx\equiv x\mod\ell^m\right\}.\]
\noindent Alors
\[G(m)\cdot\Sp_{2g}(\Z_{\ell})=\left\{M\in\GSp_{2g}((\Z_{\ell})\ |\ \lambda(M)\equiv 1\mod\ell^m\right\}.\]
\end{cor}
\demo D'apr\`es le lemme \ref{lemcong} le membre de gauche est inclus dans le membre de droite. Mais, si $M\in\GSp_{2g}(\Z_{\ell})$ et $\lambda=\lambda(M)$, la matrice
$\begin{pmatrix} I_g& 0\cr 0& \lambda(M)^{-1} I_g\cr\end{pmatrix}M$ est dans $\Sp_{2g}(\Z_{\ell})$. Si de plus $\lambda\equiv 1\mod\ell^m$, on constate que $\begin{pmatrix} I& 0\cr 0& \lambda I\cr\end{pmatrix}$ appartient \`a $G(m)$.    $\findemo$

\medskip

\begin{lemme}\label{lemmeA} Soit $G_0:=\left\{\begin{pmatrix} I& 0\cr 0& \lambda I\cr\end{pmatrix}\in\GL_{2g}(\Z_{\ell})\ |\ \lambda\in \Z_{\ell}^{\times}\right\}$, alors
$$G_0\cdot\Sp_{2g}(\Z_{\ell})=\GSp_{2g}(\Z_{\ell}).$$
\end{lemme}
\demo  Soit $M\in\GSp_{2g}(\Z_{\ell})$ de multiplicateur $\lambda(M)$. La matrice
$\begin{pmatrix} I_g& 0\cr 0& \lambda(M)^{-1} I_g\cr\end{pmatrix}M$ est dans $\Sp_{2g}(\Z_{\ell})$. $\findemo$

\medskip

\noindent Nous rassemblons maintenant quelques r\'esultats classiques sur les groupes symplectiques, leurs sous-groupes distingu\'es et leurs automorphismes.

\begin{lemme}\label{gsimple} Soit $g\geq 1$ et $K$ un corps; on exclut les cas $g=1$, $K=\F_2$, $\F_3$ ou $\F_4$ et $g=2$, $K=\F_2$, le seul sous-groupe normal non trivial de $\Sp_{2g}(K)$ est son centre $\{\pm 1\}$; les $K$-automorphismes de $\Sp_{2g}(K)$ sont tous int\'erieurs  et ceux de $\PSp_{2g}(K)$ proviennent par quotient des pr\'ec\'edents.
\end{lemme} 
\demo Voir Dieudonn\'e~\cite{dieu}, Chap. IV paragraphe 3 et Chap. IV paragraphe 6. $\findemo$

\medskip

\rem\label{automspgalg} Il est clair (c'est en fait plus facile \`a montrer) que tout automorphisme du groupe alg\'ebrique $\Sp_{2g}$ est induit par un automorphisme int\'erieur. On peut aussi \'etendre ce lemme au groupe des similitudes symplectiques.

\begin{lemme}\label{pgsimple} Soit $g\geq 1$ et $K$ un corps comme dans le lemme \ref{gsimple}. Les $K$-automorphismes de $\PGSp_{2g}(K)$ sont tous int\'erieurs.
\end{lemme}
\demo Le groupe $\PSp_{2g}(K)$ est le sous-groupe des commutateurs  de $\PGSp_{2g}(K)$, donc, si $\phi$ est un automorphisme  de $\PGSp_{2g}(K)$, sa restriction \`a $\PSp_{2g}(K)$ induit un automorphisme de $\PSp_{2g}(K)$. Or le lemme \ref{gsimple} nous indique que de tels automorphismes sont int\'erieurs et proviennent par quotient d'un automorphisme (n\'ecessairement int\'erieur) de $\Sp_{2g}(K)$. Par ailleurs, le quotient de $\PGSp_{2g}(K)$ par $\PSp_{2g}(K)$ n'est autre que le groupe $K^{\times}/K^{\times 2}$  comme on le voit en \'ecrivant le diagramme commutatif suivant~:
\[
\xymatrix{
                                                &                        1                                \ar[d]                                                                &                                                                                        1                \ar[d]                                                                &                                1                        \ar[d]                                                                                                                                                                                                &                \\
1         \ar[r]        &        \{\pm 1\}\ar[r]\ar[d]                                                        & K^{\times}\ar[r]^{x\mapsto x^2}\ar[d]        &        \left(K^{\times}\right)^{2}\ar[d]\ar[r]                                                                                & 1        \\
1         \ar[r]         & \Sp_{2g}(K)\ar[r]^{}\ar[d]        &        \GSp_{2g}(K)\ar[d]\ar[r]                                                        & K^{\times}\ar[d]\ar[r]                                                                                                                                                & 1        \\
1         \ar[r]         & \PSp_{2g}(K)\ar[r]^{}\ar[d]&        \PGSp_{2g}(K)\ar[d]\ar[r]                                                        & K^{\times}/\left(K^{\times}\right)^2\ar[d]\ar[r]        & 1        \\
                                                &                        1                                                                                                                        &                                                                                        1                                                                                                        &                                1                                                                                                                                                                                                                                                &        }
\]
\noindent Notamment tout automorphisme $\phi$ de $\PGSp_{2g}(K)$ induit un automorphisme $\phi_{|}$ sur $\PSp_{2g}(K)$. Soit $y\in \Sp_{2g}(K)$ tel que 
\[\forall x\in \PSp_{2g}(K),\ \ \phi_{|}(x)=\bar yx\bar y^{-1}.\]
\noindent D\'efinissons $\psi\in \Aut(\PGSp_{2g}(K))$ par $\psi(x)=\bar yx\bar y^{-1}.$ Montrons que $\phi\circ\psi^{-1}$ est l'identit\'e sur $\PGSp_{2g}(K)$. On sait d\'ej\`a que $(\phi\circ\psi^{-1})_{|\PSp_{2g}(K)}=id$ et on conclut gr\^ace au lemme \ref{aut1}, en observant que le centralisateur de $\PSp_{2g}(K)$ dans $\PGSp_{2g}(K)$ est trivial. 
$\findemo$

\medskip

\noindent La situation pour $\Sp_{2g}(\Z_{\ell})$ est l\'eg\`erement diff\'erente puisque les groupes de congruence
$\Gamma_n:=\Ker\{\Sp_{2g}(\Z_{\ell})\rightarrow \Sp_{2g}(\Z/\ell^n\Z)\}$ sont des sous-groupes normaux.

\begin{lemme}
Les sous-groupes normaux non triviaux de $\Sp_{2g}(\Z_{\ell})$ sont les sous-groupes $\Gamma_n$ et les sous-groupes $\tilde{\Gamma}_n$ engendr\'es par $\Gamma_n$ et $\{\pm 1\}$.
\end{lemme}

\demo Voir~\cite{kling}. $\findemo$

\begin{lemme}\label{gsimpleloc} Soit $g\geq 2$ et $R$ un anneau local de corps r\'esiduel $k$, on suppose la caract\'eristique de $k$ diff\'erente de $2,3$ ou $5$. Un automorphisme  $\Sp_{2g}(R)$ est un automorphisme int\'erieur, \'even\-tuelle\-ment multipli\'e par un caract\`ere central $\chi:\Sp_{2g}(R)\rightarrow\{\pm id\}$, c'est-\`a-dire donn\'e par $\phi(m)=\chi(m)nmn^{-1}$.
\end{lemme} 
\demo Voir  McQueen-McDonald~\cite{mcmc} pour le cas $g\geq 3$ et  Bloshchitsyn \cite{blosch} pour le cas $g=2$. $\findemo$

\begin{cor}
Soit $\ell>5$ et $R=\Z_{\ell}$ ou $R=\Z/\ell^m\Z$;  les automorphismes de   $\Sp_{2g}(R)$ sont tous des automorphismes int\'erieurs.
\end{cor}

\demo En effet un homomorphisme $\chi:\Sp_{2g}(R)\rightarrow\{\pm id\}$ se factorise dans les deux cas par un homomorphisme $\bar{\chi}:\Sp_{2g}(\F_{\ell})\rightarrow\{\pm id\}$ dont le noyau est un sous-groupe distingu\'e qui ne peut \^etre d'indice 2 et est donc \'egal au groupe $\Sp_{2g}(\F_{\ell})$ entier.
$\findemo$

\begin{lemme}\label{ssgsp} Soit $\ell>3$ premier.
Soit $H$ un sous-groupe ferm\'e de $\Sp_{2g}(\Z_{\ell})$ se projetant surjectivement sur $\Sp_{2g}(\F_{\ell})$, alors $H=\Sp_{2g}(\Z_{\ell})$ .
\end{lemme}

\demo
Voir Serre \cite{serremfv}, lemme 1 page 52.
$\findemo$

\medskip

\noindent On peut ais\'ement \'etendre ce r\'esultat \`a la situation produit : 

\begin{lemme}\label{relev} Soit $\ell\geq 5$ un nombre premier et soit $H$ un sous-groupe ferm\'e de $\GSp_{2g}(\Z_{\ell})\times\GSp_{2g}(\Z_{\ell})$ tel que son image dans $\GSp_{2g}(\F_{\ell})\times\GSp_{2g}(\F_{\ell})$ contient $\Sp_{2g}(\F_{\ell})\times\Sp_{2g}(\F_{\ell})$. Alors $H$ contient $\Sp_{2g}(\Z_{\ell})\times\Sp_{2g}(\Z_{\ell})$.
\end{lemme}

\demo Il s'agit d'une adaptation imm\'ediate du lemme 10 de \cite{serreim72} en appliquant le lemme \ref{ssgsp} du pr\'esent article en lieu et place du lemme 3 de \cite{mg}. Pour la commodit\'e du lecteur nous r\'edigeons la preuve : soit $H'$ l'adh\'erence du groupe des commutateurs de $H$. C'est un sous-groupe de $\Sp_{2g}(\Z_{\ell})\times\Sp_{2g}(\Z_{\ell})$ et son image par r\'eduction modulo $\ell$ contient le groupe d\'eriv\'e de $\Sp_{2g}(\F_{\ell})\times\Sp_{2g}(\F_{\ell})$, groupe d\'eriv\'e qui n'est autre que $\Sp_{2g}(\F_{\ell})\times\Sp_{2g}(\F_{\ell})$ tout entier car $\ell\geq 5$.  Il s'agit finalement de v\'erifier que $H'=\Sp_{2g}(\Z_{\ell})\times\Sp_{2g}(\Z_{\ell})$. Soit $X$ l'intersection de $H'$ et de $\Sp_{2g}(\Z_{\ell})\times\{1\}$, et soit $Y$ l'ensemble des \'el\'ements de $H$ dont la seconde coordonn\'ee est congrue \`a $1$ modulo $\ell$. Clairement, $X\subset Y$ et le quotient $Y/X$ est un pro-$\ell$-groupe. Soient $\tilde{Y}$ et $\tilde{X}$ les i!
 mages respectives de la premi\`ere composante de $Y$ et de $X$ dans $\Sp_{2g}(\F_{\ell})$. Par hypoth\`ese, on a $\tilde{Y}=\Sp_{2g}(\F_{\ell})$. D'autre part, $\tilde{Y}/\tilde{X}$ est isomorphe \`a un quotient de $Y/X$, donc est un $\ell$-groupe. Or $\Sp_{2g}(\F_{\ell})$ n'a pas de sous-groupe distingu\'e, autre que lui-m\^eme, d'indice une puissance de $\ell$, donc $\tilde{X}=\tilde{Y}=\Sp_{2g}(\F_{\ell})$. Le lemme \ref{ssgsp} entra\^ine que $X=\Sp_{2g}(\Z_{\ell})\times\{1\}$. Ainsi $H'$ contient le premier facteur du produit $\Sp_{2g}(\Z_{\ell})\times\Sp_{2g}(\Z_{\ell})$, et on montre de m\^eme qu'il contient le second. Donc $H'=\Sp_{2g}(\Z_{\ell})\times\Sp_{2g}(\Z_{\ell})$. $\findemo$

\medskip

\noindent Le calcul suivant combine le r\'esultat classique (sur $\F_{\ell}$) avec le lemme de Hensel (lemme \ref{hensel}).

\begin{lemme}\label{cardinal}
L'ordre des groupes  $\Sp_{2g}(\Z_{\ell})/\Gamma_n$ est donn\'e par:
$$\left| \Sp_{2g}(\Z_{\ell})/\Gamma_n\right|=\ell^{(2g^2+g)(n-1)}\left| \Sp_{2g}(\F_{\ell})\right|=\ell^{(2g^2+g)(n-1)+g^2}\prod_{i=1}^g(\ell^{2i}-1).$$
\end{lemme}

\noindent On remarquera que l'ordre d'un $\ell$-sous-groupe de Sylow est $\ell^{(2g^2+g)(n-1)+g^2}$ et l'indice d'un tel sous-groupe est $m:=\prod_{i=1}^g(\ell^{2i}-1)$ et v\'erifie $6\ell^{g(g+1)}/\pi^2\leq m\leq \ell^{g(g+1)}$. On peut en d\'eduire le corollaire suivant.

\begin{cor}\label{corcardinal} Soit $C_0>0$. Il existe une constante $C_1>0$  telle que si $ \Sp_{2g_1}(\Z_{\ell})/\Gamma_{n_1}$
et $ \Sp_{2g_2}(\Z_{\ell})/\Gamma_{n_2}$ contiennent des sous-groupes isomorphes d'indice $\leq C_0$, alors, ou bien $\ell\leq C_1$, ou bien $g_1=g_2$ et $n_1=n_2$.
\end{cor} 

\subsection{Les groupes $P_{r,s}$}

Pour d\'ecrire la dimension des stabilisateurs de sous-module du module de Tate, on introduit dans ce paragraphe les groupes alg\'ebriques qui seront, \`a conjugaison pr\`es, leurs enveloppes alg\'ebriques.
\medskip

\noindent  Dans les lemmes suivants on note $e_1,\dots,e_{2g}$ une base symplectique (\textit{i.e.} pour tout $1\leq i,j\leq g$, $e_i\cdot e_{g+i}=+1$ et $e_i\cdot e_j=0$ si $|i-j|\not= 0$).

\begin{lemme}\label{defpr} Soit $1\leq r\leq g$ et soit $P_{r}$ le sous-groupe de $\Sp_{2g}$ fixant les vecteurs $e_1,\dots,e_r$, c'est-\`a-dire
\[P_r:=\left\{M\in\Sp_{2g}\;|\; Me_i=e_i,\; i\in[1,r]\right\}.\]
\noindent Alors, $P_r$ est un sous-groupe alg\'ebrique de $\Sp_{2g}$ sur $\Z$. De plus  $P_r$ est lisse  sur $\F_{\ell}$ pour tout premier $\ell$ et de codimension
\[\codim P_r=2rg-\frac{r(r-1)}{2}.\]
\end{lemme}
\demo Il est clair que $P_r$ est un sous-groupe alg\'ebrique sur $\Z$~;~on peut donc calculer sa codimension en calculant celle de son alg\`ebre de Lie. Celle-ci est compos\'ee des matrices de $\ssp_{2g}$ dont les $r$ premi\`eres colonnes sont nulles, c'est-\`a dire de la forme
\[\begin{pmatrix}
   0_{g,r}        &        A                                & B                                        \\
  0_{r,r}                & 0_{r,g-r}        & 0_{r,g}                        \\
  0_{g-r,r}        & C                                & -^t\!\!A  
      \end{pmatrix}\]
\noindent avec $A$ matrice $g\times(g-r)$, $C$ matrice $(g-r)\times(g-r)$ sym\'etrique et $B$  matrice $g\times g$ sym\'etrique. L'\'enonc\'e en d\'ecoule ais\'ement.
$\findemo$

\medskip

\begin{lemme}\label{defprs} Soit $1\leq s\leq r\leq g$.
D\'efinissons $P_{r,s}$ le sous-groupe de $\Sp_{2g}$ fixant les vecteurs $e_1,\dots,e_r$ et les vecteurs  $e_{g+1},\dots,e_{g+s}$, c'est-\`a-dire
$$P_r:=\left\{M\in\Sp_{2g}\;|\; Me_i=e_i,\; i\in[1,r]\cup[g+1,g+s]\right\},$$
alors, $P_{r,s}$ est un sous-groupe alg\'ebrique de $\Sp_{2g}$ sur $\Z$. De plus   $P_{r,s}$ est de codimension
\[\codim P_{r,s}=2sg+2rg-rs-\frac{r(r-1)}{2}-\frac{s(s-1)}{2}.\]
De plus les groupes $P_{r,s}$ sont lisses sur $\F_{\ell}$ pour tout premier $\ell$.
\end{lemme}
\demo 
 Il est clair que $P_{r,s}$ est un sous-groupe alg\'ebrique;  on peut donc calculer sa dimension en calculant celle de son alg\`ebre de Lie. Celle-ci est compos\'ee des matrices de $\ssp_{2g}$ dont les $r$ premi\`eres colonnes sont nulles, ainsi que les colonnes $g+1,\dots,g+s$, c'est-\`a dire de la forme
\[\begin{pmatrix}
  0_{s,r}                & 0_{s,g-r}        & 0_{s,s}                & 0_{s,g-s}        \\
  0_{g-s,r}        &        A                                & 0_{g-s,s}        & B                                \\
  0_{r,r}                & 0_{r,g-r}        & 0_{r,s}                & 0_{r,g-s}        \\
  0_{g-r,r}        & C                                & 0_{g-r,s}        & -^t\!\!A  
      \end{pmatrix}\]
\noindent avec $A$ matrice $(g-s)\times(g-r)$, $C$ matrice $(g-r)\times(g-r)$ sym\'etrique et $B$  matrice $(g-s)\times (g-s)$ sym\'etrique. L'\'enonc\'e en d\'ecoule ais\'ement.
$\findemo$

\medskip

\rem Notons que le cas $s=0$ revient \`a identifier ci-dessus $P_{r,0}$ et $P_r$.

\section{Modules isotropes et propri\'et\'e $\mu$}

Comme nous supposons $A$ munie d'une polarisation principale, nous pouvons identifier l'accouplement de Weil d\'efini sur $A[\ell^n]\times A^{\vee}[\ell^n]$ (resp. sur $T_{\ell}(A)\times T_{\ell}(A^{\vee})$) \`a une application bilin\'eaire altern\'ee~: $e_{\ell^n}:A[\ell^n]\times A[\ell^n]\rightarrow\mu_{\ell^n}$ (resp. $\langle{.},{.}\rangle:T_{\ell}(A)\times T_{\ell}(A))\rightarrow T_{\ell}({\bf G}_{\rm m})$).

\subsection{Modules isotropes}

\defi Soit $M$ un sous-$\Z_{\ell}$-module de $\T_{\ell}(A)$. Nous dirons que $M$ \textit{satur\'e} si
\[\forall x\in \T_{\ell}(A),\ \ \ \ell x\in M\Rightarrow x\in M.\]

\begin{lemme}\label{sat} Soit $M$ un sous-$\Z_{\ell}$-module satur\'e de $\T_{\ell}(A)$. Notons $\pi : \T_{\ell}(A)\rightarrow \T_{\ell}(A)/\ell \T_{\ell}(A)$ le morphisme de r\'eduction modulo $\ell$. Soit $(\ee_m,\ldots,\ee_{2g})$ une base d'un suppl\'ementaire du $\F_{\ell}$-espace vectoriel $\pi(M)$. Alors tout rel\`evement  $(e_m,\ldots,e_{2g})$ de cette base, est une base d'un suppl\'ementaire de $M$ dans $\T_{\ell}(A)$.
\end{lemme}
\demo Avec les notations de l'\'enonc\'e, notons tout d'abord que la famille $(e_m,\ldots,e_{2g})$ est libre : si $\sum_i \lambda_ie_i=0$ avec $\lambda_i\in \Z_{\ell}$. Par la projection $\pi$ on en d\'eduit alors que les $\lambda_i$ sont dans $\ell\Z_{\ell}$, c'est \`a dire de la forme, $\lambda_i=\ell\mu_i$. En simplifiant par $\ell$, nous obtenons donc la relation $\sum_i \mu_ie_i=0$ et par it\'eration, on conclut que~:~ $\forall i$, $\lambda_i=0$. Notons $N$ le $\Z_{\ell}$-module engendr\'e par cette famille. On a : 
\[\pi(M+N)=\pi(M)+\pi(N)=\T_{\ell}(A)/\ell \T_{\ell}(A).\]
\noindent Or, si $L\subset \T_{\ell}(A)$ est un sous-$\Z_{\ell}$-module qui se projette surjectivement sur $\T_{\ell}(A)/\ell \T_{\ell}(A)$, alors $L=\T_{\ell}(A)$. En effet, par le th\'eor\`eme des diviseurs \'el\'ementaires, il existe une base $(f_1,\ldots,f_{2g})$ de $\T_{\ell}(A)$, un entier $2g\geq r\geq 0$ et des entiers $n_1\geq \ldots\geq n_r\geq 0$, tels que
\[L=\bigoplus_{i=1}^r\Z_{\ell}\ell^{n_i}f_i.\]
\noindent Pour que $\pi(L)$ soit un $\F_{\ell}$-espace vectoriel de dimension $2g$, on voit que, n\'ecessairement, $r=2g$ et $n_1=\ldots=n_{2g}=0$. 

\noindent Dans notre situation il nous reste donc \`a v\'erifier que la somme $M+N=\T_{\ell}(A)$ est directe. C'est ici que nous allons utiliser l'hypoth\`ese de saturation. Soit $x\in M$ et $\lambda_i\in\Z_{\ell}$ tels que $x=\sum_{i=m}^{2g}\lambda_ie_i$. On a 
\begin{align*}
x-\sum_{i=m}^{2g}\lambda_ie_i=0        & \Rightarrow \pi(x)-\sum_{i=m}^{2g}\overline{\lambda}_i\ee_i=0\\
                                                                                                                & \Rightarrow \pi(x)=0\ \text{ et }\ \forall i\geq m,\ \lambda_i\in\ell\Z_{\ell}.
\end{align*}
\noindent Donc $x\in M\cap\ell \T_{\ell}(A)=\ell M$ car $M$ est satur\'e. Finalement $x=\ell y$ avec $y\in M$ et $\lambda_i=\ell\mu_i$. En simplifiant par $\ell$ on en d\'eduit $y-\sum_{i=m}^{2g}\mu_ie_i=0$. Par it\'eration on obtient finalement $x=0$ et $\forall i\geq m$, $\lambda_i=0$.$\findemo$

\medskip

\begin{lemme}\label{isotrope} Soit $H_{\infty}$ un sous-$\Z_{\ell}$-module isotrope maximal de $\T_{\ell}(A)$. Alors $H_{\infty}$ est satur\'e.
\end{lemme}
\demo Notons $\langle \ \cdot\ \rangle$ la forme bilin\'eaire altern\'ee sur $\T_{\ell}(A)$. Soit $x\in H_{\infty}$ tel qu'il existe $y\in \T_{\ell}(A)$ tel que $x=\ell y$. Montrons que $y\in H_{\infty}$. Pour tout $z\in H_{\infty}$, on a
\[ 0=\langle x \cdot z\rangle=\langle \ell y \cdot z\rangle=\ell\langle y \cdot z\rangle.\]
\noindent Donc pour tout $z\in H_{\infty}$, on a $\langle y \cdot z\rangle=0$. Comme $H_{\infty}$ est isotrope maximal, ceci implique que $y\in H_{\infty}$.$\findemo$

\medskip

\begin{lemme}\label{h2}Soit $(e_1,\ldots,e_g)$ une base d'un sous-$\Z_{\ell}$-module isotrope maximal $H_{\infty}$ de $\T_{\ell}(A)$. Il existe un suppl\'ementaire $H_{\infty}'$ isotrope maximal et une base $(e_{g+1},\ldots,e_{2g})$ de celui-ci de sorte que
dans la d\'ecomposition $\T_{\ell}(A)=H_{\infty}\oplus H'_{\infty}$ selon la base $(e_1,\ldots,e_{2g})$, la forme symplectique s'\'ecrit comme la forme canonique $J$.
\end{lemme}
\demo Notons $\pi : \T_{\ell}(A)\rightarrow \T_{\ell}/\ell \T_{\ell}(A)$ la r\'eduction modulo $\ell$. Soit $(e_1,\ldots,e_{g})$ une base de $H_{\infty}$. Par le lemme \ref{isotrope} pr\'ec\'edent, $H_{\infty}$ est satur\'e. Nous pouvons donc appliquer le lemme \ref{sat}~: notons $(\ff_1,\ldots,\ff_g)$ une base d'un suppl\'ementaire $\overline{H'}$ totalement isotrope de $\pi(H_{\infty})$ dans $\T_{\ell}/\ell \T_{\ell}(A)\simeq \F_{\ell}^{2g}$ telle que dans la base $\{\overline{e}_1,\ldots,\overline{e}_g,\ff_1,\ldots,\ff_g\}$ de ce $\F_{\ell}$-espace-vectoriel, la matrice de la forme s'\'ecrit $J$; il nous suffit de trouver une famille $(\tf_1,\ldots,\tf_g)$ relevant les $\ff_i$, telle que 
\[\forall 1\leq i,j\leq g,\ \ \ \langle e_i\cdot \tf_j\rangle=\delta_{ij},\ \ \text{ et }\ \ \langle \tf_i\cdot \tf_j\rangle=0.\]
Par le lemme \ref{sat}, le $\Z_{\ell}$-module $H'_{\infty}$ engendr\'e par les $\tf_i$ r\'epondra au probl\`eme. Nous allons obtenir les $\tf_i$ par approximations successives modulo $\ell^n$.

\medskip

\noindent Notons $(f_1,\ldots,f_g)$ un rel\`evement des $\ff_i$ dans $\T_{\ell}(A)$ ; la famille $(f_1,\ldots,f_g)$ est une base d'un suppl\'ementaire $H'$ de $H_{\infty}$ dans $\T_{\ell}(A)$ telle que (par choix des $\ff_i$) dans la d\'ecomposition modulo $\ell$, on a $\F_{\ell}^{2g}=\pi(H_{\infty})\oplus \overline{H'}$, la forme symplectique s'\'ecrit comme la forme canonique $J$. Ceci nous donne donc une solution modulo $\ell$.

\medskip 

\noindent Voyons maintenant comment utiliser cette solution modulo $\ell$ pour obtenir une solution modulo $\ell^2$. Notons $(\tf_1,\ldots,\tf_g)\in \T_{\ell}(A)^g$ une solution potentielle modulo $\ell^2$, \textit{i.e.} telle que
\begin{equation}\label{e1}
\forall 1\leq i,j\leq g,\ \ \ \tf_i=f_i\mod\ell\ \ \text{ et }\ \ \langle e_i\cdot \tf_j\rangle=\delta_{ij}\mod \ell^2,\ \ \text{ et }\ \ \langle \tf_i\cdot \tf_j\rangle=0 \mod \ell^2.
\end{equation}
\noindent Il s'agit donc de prouver que cet ensemble d'\'equations admet une solution. On a 
\[\forall 1\leq i\leq g,\ \  \ \exists h_i\in \T_{\ell}(A),\ \ \ \tf_i=f_i+\ell h_i,\]
\noindent et on cherche donc $h_i$ solution du syst\`eme d'\'equations (\ref{e1}). Pour tout $i,j$, il existe $y_{ij}\in\Z_{\ell}$ tels que 
\[ \langle e_i\cdot f_j\rangle=\delta_{ij}+\ell y_{ij}.\]
\noindent Pour tout $i,j\in\{1,\ldots,g\}$, on a donc
\[\langle e_i\cdot \tf_j\rangle= \langle e_i\cdot f_j\rangle+\ell \langle e_i\cdot h_j\rangle=\delta_{ij}+\ell \left(y_{ij}+ \langle e_i\cdot h_j\rangle\right).\]
\noindent La premi\`ere partie du syst\`eme (\ref{e1}) se r\'e\'ecrit donc sous la forme
\begin{equation}\label{e2}
\forall i,j,\ \ \ y_{ij}+ \langle e_i\cdot h_j\rangle=0\mod \ell.
\end{equation}
\noindent Par ailleurs, 
\[\forall 1\leq i,j\leq g,\ \ \exists \alpha_{ij}\in\Z_{\ell},\ \ \ \langle f_i\cdot f_j\rangle=\ell\alpha_{ij},\]
\noindent et la matrice $(\alpha_{ij})_{i,j}$ est antisym\'etrique. Donc, on a
\begin{align*}
 \langle \tf_i\cdot \tf_j\rangle=\langle f_i+\ell h_i\cdot f_j+\ell h_j\rangle        & = \langle f_i\cdot f_j\rangle+\ell \left(\langle f_i\cdot h_j\rangle+\langle h_i\cdot f_j\rangle\right)\mod \ell^2\\
                                                                                                                                                                                                                                                                                                                         & = \ell \left(\alpha_{ij}+\langle f_i\cdot h_j\rangle+\langle h_i\cdot f_j\rangle\right)\mod \ell^2
\end{align*}
\noindent Ainsi, la seconde partie du syst\`eme (\ref{e1}) se r\'e\'ecrit
\begin{equation}\label{e3}
\forall 1\leq i<j\leq g,\ \ \ \alpha_{ij}+\langle f_i\cdot h_j\rangle+\langle h_i\cdot f_j\rangle=0\mod \ell.
\end{equation}
\noindent \'Ecrivons les inconnues $h_i$ dans la base $(e_1,\ldots,e_g,f_1,\ldots,f_g)$ : 
\[\forall 1\leq i\leq g,\ \ \ h_i=\sum_{k=1}^g h_i^ke_k+\sum_{k=1}^gh_i^{g+k}f_k,\]
\noindent avec les $h_i^k$ dans $\Z_{\ell}$. Avec ces notations, et en utilisant que dans la base $(\ee_1,\ldots,\ee_g,\ff_1,\ldots\ff_g)$ la matrice de $\langle \ \cdot\ \rangle$ est la matrice $J$, on a
\begin{equation}\label{e4}
(\ref{e2})\iff\forall 1\leq i,j\leq g,\ \ \ h_j^{g+i}=-y_{ij}\mod\ell,
\end{equation}
\noindent et
\begin{equation}\label{e5}
(\ref{e3})\iff\forall 1\leq i<j\leq g,\ \ \ \alpha_{ij}-h_j^i+h_i^j=0\mod\ell.
\end{equation}
\noindent Le syst\`eme (\ref{e4}) d\'etermine de mani\`ere unique modulo $\ell$ les composantes $h_i^{g+k}$ pour tout $i,k\in\{1,\ldots,g\}$. Le syst\`eme (\ref{e5}) ne fait intervenir que les composantes $h_i^j$ avec $j\leq g$ et peut se r\'e\'ecrire sous la forme
\[(\ref{e5})\iff \begin{cases}
\forall j\geq 2,\ \ h_1^j=\alpha_{1j}+h_j^1\mod \ell\\
\vdots\\
h_{g-1}^g=\alpha_{g-1,g}+h_g^{g-1}\mod\ell.
\end{cases}\]
\noindent Sous cette derni\`ere forme on voit imm\'ediatement qu'il existe des solutions.

\medskip

\noindent Le m\^eme calcul, en rempla\c{c}ant $\ell$ par $\ell^n$ et $\ell^{2}$ par $\ell^{n+1}$ au d\'epart, montre qu'\'etant donn\'ee une solution modulo $\ell^n$, on en d\'eduit une solution modulo $\ell^{n+1}$ qui est compatible (\textit{i.e.} se r\'eduit modulo $\ell^n$ en la solution modulo $\ell^n$ dont on est parti). Ceci nous assure donc de l'existence d'une solution $(\tf_1,\ldots,\tf_g)\in \T_{\ell}(A)^g$ comme recherch\'ee. $\findemo$

\medskip

\begin{lemme}\label{h1} Soit $H$ un sous-espace, totalement isotrope, de $A[\ell]$. Notons $\pi : \T_{\ell}(A)\rightarrow A[\ell]$ la projection canonique. Il existe un sous-$\Z_{\ell}$-module, $H_{\infty}$ de $\T_{\ell}(A)$, totalement isotrope, tel que $\pi(H_{\infty})=H$.
\end{lemme}
\demo Comme au lemme pr\'ec\'edent on raisonne par approximations successives. Con\-si\-d\'e\-rons $(\ee_1,\ldots,\ee_r)$ une base du $\F_{\ell}$-espace vectoriel $H$ que l'on compl\`ete en une base de $\T_{\ell}(A)/\ell \T_{\ell}(A)$ : $(\ee_1,\ldots,\ee_g,\ff_1,\ldots,\ff_g)$ telle que la forme $\langle\ \cdot\ \rangle$ ait pour matrice $J$ dans cette base. On commence par remonter $H$ modulo $\ell^2$ : on cherche des vecteurs $\te_1,\ldots,\te_r\in \T_{\ell}(A)$ tels que $\te_i=\ee_i\mod\ell$ et tels que 
\begin{equation}\label{ee1}
\forall 1\leq i,j\leq r,\ \ \ \langle \te_i\cdot\te_j\rangle=0\mod\ell^2.
\end{equation}
\noindent On fixe donc $e_1,\ldots,e_g,f_1,\ldots,f_g$ des rel\`evements quelconques des $\ee_i,\ff_j$ dans $\T_{\ell}(A)$, et on pose, pour tout $i\leq r$,  $\te_i=e_i+\ell h_i$. Notons que l'on a $\langle e_i\cdot e_j\rangle=\ell\alpha_{ij}$ o\`u la matrice $\alpha_{ij}$ est antisym\'etrique. Ainsi, on a
\[\langle \te_i\cdot\te_j\rangle=\ell\left(\alpha_{ij}+\langle e_i\cdot h_j\rangle+\langle h_i\cdot e_j\rangle\right)\mod\ell^2.\]
\noindent Le syst\`eme (\ref{ee1}) est donc \'equivalent \`a
\[\forall 1\leq i,j\leq r,\ \ \alpha_{ij}+\langle e_i\cdot h_j\rangle+\langle h_i\cdot e_j\rangle=0\mod\ell.\]
\noindent Cec syst\`eme est un sous-syst\`eme du syst\`eme (\ref{e3}) du lemme \ref{h2} pr\'ec\'edent. En particulier on peut trouver une solution. On conclut alors comme au lemme pr\'ec\'edent.$\findemo$

\medskip

\defi Soit $H\subset A[\ell^{\infty}]$ un sous-groupe fini. Nous dirons que \textit{$H$ est totalement isotrope} si pour tout points $P,Q$ de $H$, de m\^eme ordre $\ell^n$, on a
\[e_{\ell^n}(P,Q)=1,\]
\noindent o\`u $e_{\ell^n}$ d\'esigne l'accouplement de Weil sur $A[\ell^n]$.

\medskip

\noindent Notons que si $H$ est totalement isotrope au sens pr\'ec\'edent, alors son sous-groupe des points de $\ell$-torsion est totalement isotrope dans le $\F_{\ell}$-espace vectoriel $A[\ell]$. De plus, le lemme pr\'ec\'edent se g\'en\'eralise imm\'ediatement \`a un tel groupe $H$.

\medskip

\begin{lemme}\label{hh2} Soit $H\subset A[\ell^{\infty}]$ un sous-groupe fini, totalement isotrope. Notons $\ell^n$ l'exposant de $H$. Notons $\pi_n : \T_{\ell}(A)\rightarrow A[\ell^n]$ la projection canonique. Il existe un sous-groupe totalement isotrope $H_{\text{ti}}$ de $A[\ell^n]$, contenant $H$ et il existe un sous-$\Z_{\ell}$-module, $H_{\infty}$ de $\T_{\ell}(A)$, totalement isotrope, tel que $\pi_n(H_{\infty})=H_{ti}$.
\end{lemme}
\demo Le groupe $H$ est isomorphe \`a $\prod_{i=1}^r(\Z/\ell^{m_i}\Z)^{a_i}$ avec $m_1>\ldots>m_r$ et les $a_i\geq 1$. En suivant la proc\'edure indiqu\'ee dans la preuve du lemme \ref{h1} pr\'ec\'edent, on trouve un rel\`evement $H_1$ modulo $\ell^{m_{r-1}}$ de $(\Z/\ell^{m_r})^{a_r}$ tel que $\prod_{i=1}^{r-1}(\Z/\ell^{m_i}\Z)^{a_i}\times H_1$ est encore totalement isotrope. On rel\`eve ensuite la partie $(\Z/\ell^{m_{r-1}})^{a_{r-1}}\times H_{1}$ en un $H_2$ modulo $\ell^{m_{r-2}}$ tel que $\prod_{i=1}^{r-2}(\Z/\ell^{m_i}\Z)^{a_i}\times H_2$ est totalement isotrope. Par it\'eration on obtient un groupe $H_{ti}:=H_r$ totalement isotrope de la forme $(\Z/\ell^{m_1}\Z)^{\sum_{i=1}^ra_i}$. On le rel\`eve, toujours par la m\^eme proc\'edure en un $H_{\infty}$ qui convient.$\findemo$

\subsection{Propri\'et\'e $\mu$}

\noindent \'Etant donn\'e un sous-groupe $H$ fini de $A[\ell^{\infty}]$, nous introduisons \`a pr\'esent l'invariant suivant~:
\[m_1(H)=\max\left\{k\in\N\ |\ \exists n\geq 0, \ \exists P,Q \in H\text{ d'ordre }\ell^n,\ \ e_{\ell^n}(P,Q) \text{ engendre }\mu_{\ell^k}\right\}.\]
\noindent Dire que $H$ est totalement isotrope \'equivaut \`a dire que $m_1(H)=1$. De plus on peut noter que, sur la d\'efinition, il est \'evident que $m_1(H)$ est sup\'erieur \`a la valeur $m$ suivante :
\[m(H):=\max\left\{k\in\N\ |\ \exists P,Q \in H\text{ d'ordre }\ell^k,\ \ e_{\ell^k}(P,Q) \text{ engendre }\mu_{\ell^k}\right\}.\]
\noindent Lorsque $H$ est de la forme $A[\ell^n]$ alors $m_1(H)=m(H)=n$.

\medskip

\noindent  Dans le cas g\'en\'eral, si $H$ contient deux points d'ordre $\ell^n$, tel que l'accouplement de Weil de ces deux points est une racine primitive $k$-i\`eme de l'unit\'e, alors comme cet accouplement est Galois-\'equivariant, on obtient que
\[K(\mu_{\ell^k})\subset K(H).\] 

\defi Nous appelons \textit{propri\'et\'e ($\mu$)} pour une vari\'et\'e ab\'elienne le fait d'avoir, pour tout sous-groupe fini $H\subset A\left[\ell^{\infty}\right]$, l'\'egalit\'e (\`a indice fini pr\`es, born\'e uniform\'ement)~:
\[ K(\mu_{\ell^{m_1(H)}})= K(H)\cap K(\mu_{\ell^{\infty}}).\]

\medskip

\noindent Nous allons montrer que les vari\'et\'es ab\'eliennes de type $\GSp$ ont la propri\'et\'e ($\mu$). Notons que, vu notre d\'efinition de l'invariant $m_1(H)$, on a
\[K(\mu_{\ell^{m_1(H)}})\subset K(H)\cap K(\mu_{\ell^{\infty}}).\]

\begin{prop}\label{pmu} Soit $A$ une vari\'et\'e ab\'elienne de dimension $g$, d\'efinie sur un corps de nombres $K$. En notant $\delta(H):=\left(\Z_{\ell}^{\times}:\lambda(G_0(H))\right)$ o\`u $G_0(H)=\Gal\left( K(A\left[{\ell^{\infty}}\right])/K(H)\right)$, on a, \`a indice fini pr\`es, pour tout $H$ sous-groupe fini de $A[\ell^{\infty}]$,
\[ [K(H)\cap K(\mu_{\ell^{\infty}}):K]= \delta(H).\]
Si on suppose de plus que la vari\'et\'e ab\'elienne $A$ est telle que $G_{\ell}$ s'identifie, \`a indice fini pr\`es, avec $\GSp_{2g}(\Z_{\ell})$, alors, pour tout $H$ sous-groupe fini de $A[\ell^{\infty}]$, on a l'\'egalit\'e \`a indice fini pr\`es~:
\[K(H)\cap K(\mu_{\ell^{\infty}})=K(\mu_{\ell^{m_1(H)}}).\]
\end{prop}
\demo
\noindent 
Identifions le groupe de Galois $\Gal\left( K(A\left[{\ell^{\infty}}\right])/K\right)$ \`a un sous-groupe de  $\GSp(\Z_{\ell})$. Le groupe de Galois $\Gal\left( K(A\left[{\ell^{\infty}}\right])/K(\mu_{\ell^{\infty}})\right)$ s'identifie alors \`a avec $SG_{\ell}:=G_{\ell}\cap \Ker(\lambda)$. Alors  $K(H)\cap K(\mu_{\ell^{\infty}})$ est la sous-extension fix\'ee par le groupe $U$ engendr\'e par $SG_{\ell}$ et $G_0(H)$. On voit imm\'ediatement que le noyau de
$G_{\ell}\buildrel{\lambda}\over{\rightarrow}\Z_{\ell}^{\times}\rightarrow\Z_{\ell}^{\times}/\lambda(G_0(H))$ est le groupe $U$ d'o\`u le premier \'enonc\'e.
\medskip

\noindent Passons maintenant au cas d'une vari\'et\'e de type Gsp. Commen\c{c}ons par consid\'erer $H_{\infty}$ un sous-groupe isotrope maximal de $\T_{\ell}(A)$. Par le lemme \ref{h2}, on peut supposer que dans une d\'ecomposition $\T_{\ell}(A)=H_{\infty}\oplus H'_{\infty}$ la forme symplectique s'\'ecrit comme la forme canonique $J$. On voit alors ais\'ement que
\begin{align*}
\Gal\left(K(A[\ell^{\infty}])/K(H_{\infty})\right)        & =\left\{M=\begin{pmatrix} I& *\cr 0& *\cr\end{pmatrix}\in\GSp_{2g}(\Z_{\ell})\right\}\\
                                                                                                                                                                                                                & =\left\{M=\begin{pmatrix} I& S\cr 0& \lambda I\cr\end{pmatrix}\;|\; \lambda\in\Z_{\ell}^{\times}\;\text{et}\; S\;\text{sym\'etrique}\right\}
\end{align*}
\noindent D'apr\`es le lemme \ref{lemmeA}, le groupe engendr\'e par ce dernier groupe et par le groupe
$\Sp_{2g}(\Z_{\ell})=\Gal\left(K(A[\ell^{\infty}])/K(\mu_{\ell^{\infty}})\right)$ est $\GSp_{2g}(\Z_{\ell})$ tout entier. Ainsi $K(H_{\infty})\cap K(\mu_{\ell^{\infty}})=K$. Si $H$ est un sous-groupe fini de $A[\ell^{\infty}]$ totalement isotrope, dans ce cas, le lemme \ref{hh2} et ce qui pr\'ec\`ede nous permettent de conclure :  \`a indice fini pr\`es, on a $K(H)\cap K(\mu_{\ell^{\infty}})=K$.

\medskip

\noindent Soit maintenant $H$ un sous-groupe fini non-isotrope d'exposant $\ell^{r_H}$ de $A[\ell^{\infty}]$. On a 
\[ [\ell^{m_1(H)}](H) \text{ est totalement isotrope.}\]
\noindent En effet si $P$ et $Q$ sont deux points d'ordre $\ell^n$ dans $H$, alors 
\[e_{\ell^{n-m_1(H)}}(\ell^{m_1(H)}P,\ell^{m_1(H)}Q)=e_{\ell^{n}}(P,Q)^{\ell^{m_1(H)}}=1 \text{ par d\'efinition de }m_1(H).\]
\noindent En appliquant le lemme \ref{hh2}, on trouve donc un sous-groupe $H'$ contenant $[\ell^{m_1}](H)$ de m\^eme exposant et il existe un sous-$\Z_{\ell}$-module $H_{\infty}$ totalement isotrope de $\T_{\ell}(A)$ tel que,  si, pour tout entier $n\geq 1$, $\pi_n : \T_{\ell}(A)\rightarrow \T_{\ell}(A)/\ell^n\T_{\ell}(A)=A[\ell^n]$ d\'esigne la projection canonique, on a
\[\pi_{r_H}(H_{\infty})=H'.\]
\noindent Par le lemme \ref{h2}, on peut supposer que dans une d\'ecomposition $\T_{\ell}(A)=H_{\infty}\oplus H'_{\infty}$ la forme symplectique s'\'ecrit comme la forme canonique $J$. Pour tout $n\geq 1$, notons 
\[H_n:=\pi_n(H_{\infty})=H_{\infty}/H_{\infty}\cap\ell^n \T_{\ell}(A).\]
\noindent On a pour tout $n\geq 1$, $[\ell] H_{n+1}=H_n$. On peut donc poser
\[H^{\infty}=\bigcup_{n\geq 1}H_n\subset A[\ell^{\infty}].\]
\noindent De plus, on voit que le groupe de Galois correspondant \`a $H_{\infty}$ est le m\^eme que celui correspondant \`a $H^{\infty}$. On a 
\[H\subset[\ell^{m_1(H)}]^{-1}(H')=[\ell^{m_1(H)}]^{-1}(H_{r_H})\subset[\ell^{m_1(H)}]^{-1}(H^{\infty}).\]
\noindent En consid\'erant la multiplication par $\ell^{m_1(H)}$ sur $H^{\infty}$, on en d\'eduit (car $H^{\infty}$ est $\ell$-divisible)   que 
\[H\subset H^{\infty}+\ker[\ell^{m_1(H)}]=:\widetilde{H^{\infty}}.\]
\noindent Ainsi comme dans le cas totalement isotrope, on se ram\`ene \`a une situation o\`u un lemme de groupe permet de conclure : 
\[\Gal\left(K(A[\ell^{\infty}])/K(\widetilde{H^{\infty}})\right)=\left\{M\in\GSp_{2g}(\Z_{\ell})\;|\;\forall i\leq g,\ Me_{g+i}=e_{g+i}\text{ mod }\ell^{m_1(H)},\text{ et, }Me_i=e_i\right\}.\]
\noindent La m\^eme preuve que celle du corollaire \ref{corB} donne alors le r\'esultat : 
\[\Gal\left(K(A[\ell^{\infty}])/K(\widetilde{H^{\infty}})\right)\cdot \Sp_{2g}(\Z_{\ell})=\left\{M\in\GSp_{2g}((\Z_{\ell})\ |\ \lambda(M)\equiv 1\mod\ell^{m_1(H)}\right\}.\]
\noindent Notamment, on a, \`a indice fini born\'e pr\`es,
\[K(H)\cap K(\mu_{\ell^{\infty}})\subset K(\widetilde{H^{\infty}})\cap K(\mu_{\ell^{\infty}})\subset K(\mu_{\ell^{m_1(H)}}).\]
\noindent Par la remarque pr\'ec\'edant la proposition on sait d\'eja que l'autre inclusion est vraie. $\findemo$

\section{Calcul de l'invariant $\gamma(A)$ pour $A$ simple de type GSp\label{simple}}

Nous d\'emontrons maintenant un r\'esultat qui, conjugu\'e au th\'eor\`eme \ref{th1} entra\^{\i}ne le th\'eor\`eme \ref{th2}.

\begin{theo}\label{th3}  Si $A/K$ est une vari\'et\'e ab\'elienne de dimension $g$, d\'efinie sur un corps de nombres $K$, telle que~:  pour tout premier $\ell$, le groupe $G_{\ell}$ s'identifie \`a un sous-groupe de $\GSp_{2g}(\Z_{\ell})$ d'indice fini, born\'e ind\'ependamment de $\ell$ ; alors  
\[\gamma(A)=\frac{2\dim A}{\dim\MT(A)}=\frac{2g}{2g^2+g+1}.\]
\end{theo}
\demo Nous commen\c{c}ons l'argument dans le cas plus g\'en\'eral d'une vari\'et\'e ab\'elienne quelconque v\'erifiant la conjecture de Mumford-Tate forte (la d\'efinition des groupes de Hodge et de Mumford-Tate est rappel\'ee au paragraphe \ref{MTH}). Nous sp\'ecialiserons un peu plus tard la preuve au cadre $\MT(A)=\GSp_{2g}$. Pour simplifier les notations, nous supposerons ici que $\Gal(K(A[\ell^{\infty}])/K)$ s'identifie avec $\MT(\Z_{\ell})$. On commence par se ramener au cas $\ell$-adique (cf. \cite{hindry-ratazzi1} proposition 4.1)~: 
\begin{prop}\label{ladique}Soit $\alpha>0$. Pour d\'emontrer que $\gamma(A)\leq\alpha$, il suffit de montrer que : il existe une cons\-tan\-te strictement positive $C(A/K)$ ne d\'ependant que de $A/K$ telle que pour tout nombre premier $\ell$, pour tout sous-groupe fini $H_{\ell}$ de $A[\ell^{\infty}]$, on a
\begin{equation}\label{equa}
\text{Card}\left(H_{\ell}\right)\leq C(A/K)[K(H_{\ell}):K]^{\alpha}.
\end{equation}
\end{prop}

\medskip

\noindent Soit donc $H$ sous-groupe fini de $A[\ell^{\infty}]$, on pose
$$G_0(H):=\left\{M\in\MT(\Z_{\ell})\;|\;\forall x\in H,\;Mx=x\right\}.$$
et $G(H):=G_0(H)\cap\Hdg(\Z_{\ell})$. Comme groupe abstrait, $H$ est de la forme $H\simeq\prod_{i=1}^{2g}\Z/\ell^{m_i}\Z$. Notons $e_1, \ldots, e_{2g}$ un syst\`eme de g\'en\'erateurs ; les $e_i$ \'etant d'ordre respectif $\ell^{m_i}$. Notons de plus $\{\hat{e_1},\ldots,\hat{e}_{2g}\}$ une base de $\T_{\ell}(A)$ relevant la famille $\{e_i\}$, \textit{i.e.} telle que $e_i=\hat{e}_i\mod \ell^{m_i}$ pour tout $i$. On a
\[G(H)=\left\{M\in\Hdg(\Z_{\ell})\;|\; M\hat{e}_i= \hat{e}_i\mod\ell^{m_i},\ 1\leq i\leq 2g\right\}.\]

\medskip

\begin{lemme} \label{calculH1}  Soit $H$ un sous-groupe fini de $A\left[\ell^{\infty}\right]$. Notons
$\delta(H):=\left(\Z_{\ell}^{\times}:\lambda(G_0(H))\right)$. 
On a alors~:
\[[K(H):K]=(\MT(\Z_{\ell}):G_0(H))=\delta(H)(\Hdg(\Z_{\ell}):G(H)).\]
\end{lemme}
\demo La premi\`ere \'egalit\'e est donn\'ee par la th\'eorie de Galois car on a suppos\'e que $\Gal(K(A[\ell^{\infty}])/K)$ s'identifie avec $\MT(\Z_{\ell})$. La seconde \'egalit\'e est une chasse au diagramme facile.
$\findemo$

\medskip

\noindent Nous supposerons d\'esormais que $A$ est de type $\GSp$, et donc que $\Hdg(A)=\Sp_{2g}$. Quitte \`a renum\'eroter on peut supposer que les exposants $m_i$ (correspondants au $e_i$) sont ordonn\'es dans l'ordre d\'ecroissant : $m_1\geq\ldots\geq m_{2g}$. On pose alors
\[m^1:=\max \{m_i\ |\ m_i\not=0\}\ \text{ et par r\'ecurrence}\ m^{r+1}=\max\{m_i\ |\ m_i<m^r\}.\]
\noindent On obtient ainsi une suite strictement d\'ecroissante $m^1>\ldots > m^t\geq 1$ (avec $t\leq 2g$). Le groupe $H$ est isomorphe \`a $\prod_{i=1}^t\left(\Z/\ell^{m^i}\Z\right)^{a_i}$. On d\'efinit ensuite pour tout $1\leq r\leq t$, les sous-ensembles 
\[ I_{r}=\{i\in\{1,\ldots,2g\}\ |\ m_i\geq m^{r}\}\ \ \text{ de cardinal }\ \ \left|I_{r}\right|=\sum_{i=1}^{r}a_i.\]
\noindent Introduisons maintenant la suite croissante de groupes alg\'ebriques sur $\Z_{\ell}$ suivants :
\[\forall 1\leq r\leq t\ \  G_{r}:=\left\{M\in \Sp_{2g}\ |\ M\hat{e}_i=\hat{e}_i\ \ \forall i\in I_{t+1-r}\right\}.\]
\noindent On voit que
\[G(H)=\left\{M\in \Sp_{2g}(\Z_{\ell})\ |\ \forall 1\leq r\leq t\ \ M\in G_r\mod\ell^{m^{t+1-r}}\right\}.\]
\noindent Par changement de base symplectique sur $\F_{\ell}$, chacun des $G_i$ est conjugu\'e sur $\F_{\ell}$ \`a l'un des groupes $P_{r,s}$ introduits au paragraphe \ref{prs}. En posant $G=\Sp_{2g}$, on voit que, avec les notations du lemme \ref{cle}, on a 
\[G(H)=H(m^1,\ldots,m^t).\]
\noindent On va donc pouvoir appliquer le lemme \ref{cle}. Pour obtenir la valeur exacte de $\gamma(A)$ et pas seulement une majoration, nous allons traiter le cas $H$ d'exposant $\ell$ puis le cas g\'en\'eral.

\subsection{Si $H$ est d'exposant $\ell$}

\noindent Commen\c{c}ons par calculer $\delta(H)$. Le groupe $H$ est inclus dans le $\F_{\ell}$-espace vectoriel $A[\ell]$. 

\begin{lemme}Si $H$ est inclus dans un sous-espace totalement isotrope alors $\delta(H)=1$. Sinon $\delta(H)\gg\ll \ell$.
\end{lemme}
\demo Si $H$ est totalement isotrope, c'est la proposition \ref{pmu}. Si $H$ est n'est pas inclus dans un lagrangien, alors il existe $P, Q\in H$ tels que l'accouplement de Weil, $e_{\ell}(P,Q)$, soit une racine $\ell$-i\`eme primitive de l'unit\'e. Le formalisme de l'accouplement de Weil nous donne alors
\[K(\mu_{\ell})\subset K(H)\cap K(\mu_{\ell^{\infty}}).\] 
\noindent De plus on a $K(H)\subset K(A[\ell])$ et on sait par la proposition 6.8 de \cite{hindry-ratazzi1} que 
\[\left[K(A[\ell])\cap K(\mu_{\ell^{\infty}}): K(\mu_{\ell})\right]=O(1).\]
\noindent On en d\'eduit donc le r\'esultat.$\findemo$

\medskip

\noindent  Il reste maintenant \`a calculer $(\Sp(\Z_{\ell}):G(H))$. Dans notre cadre d'exposant $\ell$,  l'entier $t$ pr\'ec\'edent vaut n\'ecessairement $1$. Le lemme \ref{cle} (applicable d'apr\`es le lemme  \ref{defprs}) nous donne donc : 
\[ [K(H):K]\gg\ll \delta(H)\ell^{\codim P_{r,s}},\]
\noindent o\`u $(r,s)$ (avec \'eventuellement $s=0$) est le couple correspondant \`a $H$. 

\begin{enumerate}
\item Si $H$ est inclus dans un lagrangien, on obtient donc  l'in\'egalit\'e
\[|H|=\ell^r\ll [K(H):K]^{\gamma}\]
\noindent si et seulement si 
\[\gamma\geq \frac{r}{\codim P_r}.\]
\noindent Il est clair que lorsque $H$ varie, tous les groupes $P_r$ interviennent, donc pour les groupes totalement isotropes, l'exposant
\[\max_{1\leq r\leq g}\frac{r}{\codim P_r}\]
\noindent est admissible.
\noindent Dans notre cas, on a $\Hdg(A)=\Sp_{2g}$, donc
\[\frac{1}{2g-\frac{r}{2}+\frac{1}{2}}=\frac{r}{2rg-\frac{r(r-1)}{2}}\leq \gamma.\]
\noindent Un simple calcul montre que le membre de gauche, vu comme fonction de $r\in[1,g]$, est croissant et a pour valeur maximale
$\frac{2g}{3g^2+g}$ (pour $r=g$) qui est bien plus petit que $\gamma:= \frac{2g}{2g^2+g+1}$.

\item Si $H$ n'est pas inclus dans un lagrangien, on obtient de m\^eme l'in\'egalit\'e
\[|H|=\ell^{r+s}\ll [K(H):K]^{\gamma}\gg\ll \ell^{\gamma(1+\codim P_{r,s})}\]
\noindent si et seulement si 
\[\gamma\geq \frac{r+s}{1+\codim P_{r,s}}.\]
\noindent Finalement dans ce cas, l'exposant
\[\max_{1\leq r,s\leq g}\frac{r+s}{1+\codim P_{r,s}}\]
\noindent est admissible.
\noindent Comme $\Hdg(A)=\Sp_{2g}$, on a donc
\[[K(H):K]\gg\ll \ell^{1+\codim P_{r,s}}=\ell^{1+2sg+2rg-rs-\frac{r(r-1)}{2}-\frac{s(s-1)}{2}},\]    
\noindent et on voit que l'on a l'in\'egalit\'e
\[|H|=\ell^{r+s}\ll [K(H):K]^{\gamma}\gg\ll \ell^{\gamma(1+\codim P_{r,s})}\]
\noindent si et seulement si 
\[\frac{r+s}{1+2sg+2rg-rs-\frac{r(r-1)}{2}-\frac{s(s-1)}{2}}\leq \gamma.\]
\noindent Un simple calcul montre que le membre de gauche, vu comme fonction de $r,s\in[1,g]$, est croissante par rapport aux deux variables, et a une valeur maximale
 (pour $r=s=g$) qui est \'egale $\gamma:= \frac{2g}{2g^2+g+1}$.
\end{enumerate}

\medskip

\noindent Finalement dans le cas o\`u le groupe $H$ est d'exposant $\ell$, on voit que 
\[\alpha(A)=\frac{2g}{2g^2+g+1}\]
est un exposant admissible et de plus, c'est \'egalement la valeur minimale pour $\gamma(A)$.

\subsection{Si $H$ est quelconque}
\noindent On suppose maintenant que $H$ est d'exposant $\ell^{n}$. Il suffit, pour conclure, de montrer que la valeur $\alpha(A)$ pr\'ec\'edente est toujours admissible dans ce cas. Comme pr\'ec\'edemment, on peut appliquer le lemme \ref{cle}~:
\[\left[\Sp_{2g}(\Z_{\ell}):G(H)\right]\gg \ell^{\sum_{i=1}^td_i(m^{t+1-i}-m^{t+1-(i-1)})},\]
\noindent o\`u l'on a pos\'e $m^{t+1}=0$ et o\`u $d_i$ est la codimension de $G_i$. Les groupes alg\'ebriques $G_i$ \'etant conjugu\'es sur $\F_{\ell}$ aux $P_{r,s}$ (avec \'eventuellement $s=0$), $d_i$ est \'egalement la codimension du groupe $P_{r_i,s_i}$ correspondant. Par ailleurs, la suite des $(G_i)_i$ \'etant croissante, la suite des $(P_{r_i,s_i})_i$ l'est \'egalement. Ceci se traduit par
\[\forall i,\ \ \ r_i\geq r_{i+1}\ \ \ \text{ et }\ \ \ s_i\geq s_{i+1}.\]
\noindent Il nous reste \`a calculer la valeur de $\delta(H)$ (ou plutot une minoration de $\delta(H)$). Soit $m\geq 1$ un entier maximal tel que $H$ contient deux points $P, Q$ d'ordre $\ell^m$ tels que $e_{\ell}(\ell^{m-1}P,\ell^{m-1}Q)$ est une racine primitive $\ell$-i\`eme de l'unit\'e. Si un tel $m$ n'existe pas, posons $m:=0$. Soit ensuite $h\in\{1,\ldots,t\}$ minimal tel que $m^h\leq m$. Si $m=0$, prenons $h=t+1$. L'accouplement de Weil nous indique que 
\[\delta(H)\gg\phi(\ell^m)\geq\phi(\ell^{m^h})\gg\ll\ell^{m^{h}}.\]
\noindent Par ailleurs, si $h=t+1$ tous les $P_{r_i,s_i}$ sont tels que $s_i=0$. Si par contre $h\leq t$, alors le groupe $P_{r_h,s_h}$ est tel que $s_h\not=0$ ; pour $k\geq h$, on a $s_{t+1-k}\not=0$, et pour $k<h$, on a $s_{t+1-k}=0$. Autrement dit on a la suite d'inclusions
\[P_{r_1,s_1}\subset\ldots\subset P_{r_{t+1-h},s_{t+1-h}}\subset P_{r_{t+1-(h-1)}}\subset \ldots\subset P_{r_t}.\]
\noindent Posons
\[\delta_1=\ldots=\delta_{t+1-h}=1,\ \text{ et }\ \delta_{t+1-(h-1)}=\ldots=\delta_t=0.\]
\noindent On voit que $m^h=\sum_{i=1}^t(m^{t+1-i}-m^{t+1-(i-1)})\delta_i.$ Nous obtenons ainsi la minoration
\[[K(H):K]\gg\ell^{\sum_{i=1}^t(\m^{t+1-i}-m^{t+1-(i-1)})(\delta_i+\codim P_{r_i,s_i})}.\]
\noindent De plus, pour tout entier $k\in\{1,\ldots t\}$,
\[r_{t+1-k}+s_{t+1-k}=|I_k|=\sum_{i=1}^ka_i.\]
\noindent On aura donc $|H|=\ell^{a_1m^1+\dots+a_tm^t}\ll [K(H):K]^{\gamma}$ si \[\gamma\geq\max\left\{\frac{a_1m^1+\dots+a_tm^t}{\sum_{j=1}^t(m^j-m^{j-1})(\delta_{t+1-j}+\codim P_{r_{t+1-j},s_{t+1-j}})}\right\},\]
le maximum \'etant pris sur $m^1>\dots> m^t$ et  $0\leq s_j\leq r_j\leq g$ et $r_t+s_t=a_1+\dots+a_t\leq 2g$. Rappelons la notation $d_j=\codim P_{r_j,s_j}$. Le maximum cherch\'e est donc :
\[M=\max_{m^1\geq \dots\geq m^t}\left\{\frac{\sum_{i=1}^ta_im^i}{\sum_{i=1}^tm^i(\delta_{t+1-i}-\delta_{t+1-(i-1)}+d_{t+1-i}-d_{t+1-(i-1)})}\right\}.\]
\noindent D'apr\`es le lemme \ref{combielem} on a
\[ M=\max_{1\leq k\leq t}\left\{\frac{\sum_{i=1}^ka_i}{\sum_{i=1}^k(\delta_{t+1-i}-\delta_{t+1-(i-1)}+d_{t+1-i}-d_{t+1-(i-1)}) }\right\}\]
\noindent Soit en simplifiant : 
\[M=\max_{1\leq k\leq t}\frac{r_{t+1-k}+s_{t+1-k}}{\delta_{t+1-k}+\codim P_{r_{t+1-k},s_{t+1-k}}}.\]
Si le maximum  correspond \`a $k$ tel que $s_{t+1-k}=0$, on a alors $\delta_{t+1-k}=0$. Si le maximum  correspond \`a $k$ tel que $s_{t+1-k}>0$, on a alors $\delta_{t+1-k}=1$. On voit donc que ce maximum n'est autre que
\[\alpha(A)=\max\left\{\max_{1\leq r\leq g}\frac{r}{\codim P_r}, \max_{1\leq r,s\leq g}\frac{r+s}{1+\codim P_{r,s}}\right\}=\frac{2g}{2g^2+g+1}.\]
\noindent Ceci prouve que la valeur $\alpha(A)$ est toujours admissible et comme on a prouv\'e que pour les groupes d'exposant $\ell$, cette valeur donne \'egalement une minoration de $\gamma(A)$, on obtient le r\'esultat : pour toute vari\'et\'e ab\'elienne de type $\GSp$,  \[\gamma(A)=\alpha(A)=\frac{2g}{2g^2+g+1}.\]

\section{Mumford-Tate et Hodge pour un produit de vari\'et\'es ab\'e\-liennes de type $\GSp$\label{MTH}}

Commen\c{c}ons par rappeler la d\'efinition des groupes de Hodge et Mumford-Tate associ\'es \`a une vari\'et\'e ab\'elienne $A$ d\'efinie sur $K\subset\C$.
 On note $V=H^1(A(\C),\Q)$ le premier groupe de cohomologie singuli\`ere de la vari\'et\'e analytique complexe $A(\C)$. C'est un $\Q$-espace vectoriel de dimension $2g$. Il est naturellement muni d'une structure de Hodge de type $\{(1,0),(0,1)\}$, c'est-\`a-dire d'une d\'ecomposition sur $\C$ de $V_{\C}:=V\otimes_{\Q}\C$ donn\'ee par $V_{\C}=V^{1,0}\oplus V^{0,1}$ telle que $V^{0,1}=\overline{V^{1,0}}$ o\`u $\overline{\hspace{.1cm} \cdot \hspace{.1cm}}$ d\'esigne la conjugaison complexe. On note $\mu : \G_{m,\C}\rightarrow \GL_{V_{\C}}$ le cocaract\`ere tel que pour tout $z\in\C^{\times}$, $\mu(z)$ agit par multiplication par $z$ sur $V^{1,0}$ et agit trivialement sur $V^{0,1}$. On d\'efinit le groupe de Mumford-Tate en suivant \cite{pink1}.

\medskip

\defi Le \textit{groupe de Mumford-Tate} $\MT(A)/\Q$ de $A$ est le plus petit $\Q$-sous-groupe alg\'ebrique $G$ de $\GL_V$ (vu comme $\Q$-sch\'ema en groupes) tel que, apr\`es extension des scalaires \`a $\C$, le cocaract\`ere $\mu$ se factorise \`a travers $G_{\C}:=G\times_{\Q}\C$. Le \textit{groupe de Hodge} $\hdg(A)/\Q$ de $A$ est $(\MT(A)\cap\SL_V)^0$, la composante neutre de $\MT(A)\cap\SL_V$.

\medskip

\noindent On veut montrer le th\'eor\`eme suivant (th\'eor\`eme \ref{prophodgeprod} de l'introduction). Pour la preuve des points 2 et 3 de ce th\'eor\`eme nous adaptons au cas des vari\'et\'es ab\'eliennes de type $\GSp$, en suivant sa strat\'egie de preuve, le paragraphe 6 de \cite{serreim72}.

\medskip

\begin{theo}  {\rm \bf(= Th\'eor\`eme  \ref{prophodgeprod})} Soient $r$ et $n_1\ldots,n_r$ des entiers strictement positifs. Soient $A_i$ des vari\'et\'es ab\'e\-liennes de dimension $g_i$ non isog\`enes deux \`a deux telles que $\hdg(A_i)=\Sp_{2g_i}$. Posons $A:=A_1^{n_1}\times\dots\times A_r^{n_r}$ et, pour tout premier $\ell$, notons $\rho_{\ell,i}$ (respectivement $\rho_{\ell}=\rho_{\ell,1}\times\ldots,\rho_{\ell,r}$) les repr\'esentations $\ell$-adiques associ\'ees aux $A_i$ (respectivement \`a $A$) alors~:
\begin{enumerate}
\item l'inclusion naturelle suivante est un isomorphisme~:
$$\hdg(A)\cong\hdg\left(A_1\times\dots\times A_r\right)\hookrightarrow\Sp_{2g_1}\times\dots\times\Sp_{2g_r}.$$
\item soit $\ell$ un nombre premier. Si  pour tout $i$, on a $\rho_{\ell,i}\left(\Gal(K(A_i[\ell^{\infty}])/K(\mu_{\ell^{\infty}}))\right)\cong \Sp_{2g_i}(\Z_{\ell})$ (\`a indice fini pr\`es) 
alors, on a (\`a indice fini pr\`es)~:
$$ \rho_{\ell}\left(\Gal(K(A[\ell^{\infty}])/K(\mu_{\ell^{\infty}}))\right)\cong \Sp_{2g_1}(\Z_{\ell})\times\dots\times\Sp_{2g_r}(\Z_{\ell}).$$ 
\item si de plus l'indice fini pour chaque $A_i$ est born\'e ind\'ependamment de $\ell$, il en est de m\^eme pour $A$.
\end{enumerate}
\end{theo}

\medskip

\noindent On se ram\`ene au produit de deux facteurs gr\^ace au lemme suivant.

\begin{lemme} Soit $r\geq 1$ un entier. Soient $G_1,\dots,G_r$ des groupes alg\'ebriques semi-simples (tels que $[G_i,G_i]=G_i$ et $H$  un sous-groupe alg\'ebrique de $G_1\times\dots\times G_r$ se projetant surjectivement sur $G_i\times G_j$. Alors $H=G_1\times\dots\times G_r$.

Soient $G_1,\dots,G_r$ des groupes pro-finis tels que, pour tout sous-groupe ouvert $U\subset G_i$, l'adh\'e\-ren\-ce des commutateurs de $U$ est ouvert dans $G_i$. Soit $H$ un sous-groupe ferm\'e de $G_1\times\dots\times G_r$ se projetant surjectivement sur $G_i\times G_j$. Alors $H=G_1\times\dots\times G_r$.
\end{lemme}
\demo Voir Ribet~\cite{ribet1} lemma 3.4.$\findemo$

\medskip

\noindent Pour traiter le cas de deux facteurs on utilise le lemme classique.

\begin{lemme}\label{goursat1} \textnormal{\textbf{(Lemme de Goursat)}} Soient $H\subset G_1\times G_2$ avec $p_i(H)=G_i)$. Soit
$N_1$ (resp. $N_2$) tel que $N_1\times\{e_2\}=\Ker (p_2)\cap H$ (resp. $\{e_1\}\times N_2=\Ker (p_1)\cap H$), alors $H$ est l'image inverse du graphe dans $G_1/N_1\times G_2/N_2$ de l'isomorphisme naturel $G_1/N_1\cong G_2/N_2$.
\end{lemme}

\demo Voir Ribet~\cite{ribet1}  lemma 3.2.$\findemo$

\medskip

\noindent Nous utiliserons ci-dessous la version suivante des r\'esultats fondamentaux de Faltings \cite{falt}.

\begin{prop}\label{lemfalt} \textnormal{\textbf{(Faltings)}} Soient $A$ et $B$ deux vari\'et\'es ab\'eliennes sur un corps de nombres $K$. Notons $\rho_{\ell,A}:\Gal(\bar{K}/K)\rightarrow \Aut(\T_{\ell}(A))$, respectivement $\bar{\rho}_{\ell,A}:\Gal(\bar{K}/K)\rightarrow \Aut(A[\ell])$, la repr\'esentation associ\'ee \`a l'action de Galois sur les points de torsion de $A$. On d\'efinit de m\^eme $\rho_{\ell,B}$ et $\bar{\rho}_{\ell,B}$.
\begin{itemize}
\item Si $\rho_{\ell,A}$ est isomorphe \`a $\rho_{\ell,B}$ alors $A$ est isog\`ene \`a $B$.
\item Il existe $C_0=C_0(A,K)$ telle que si $\ell\geq C_0$ et  $\bar{\rho}_{\ell,A}\cong\bar{\rho}_{\ell,B}$ alors $A$ est isog\`ene \`a $B$.
\end{itemize}
\end{prop}
\demo Dans l'article de Faltings \cite{falt}, le premier \'enonc\'e est d\'emontr\'e dans le Korollar 2, page 361; le deuxi\`eme \'enonc\'e, pour $\bar{\rho}_{\ell}$,  peut se d\'eduire des d\'emonstrations comme cela est montr\'e par Zarhin \cite{zar1}, Corollary 5.4.5. 
$\findemo$

\medskip
\noindent{\bf Remarque.} On peut aussi d\'emontrer (voir \cite{zar2}, paragraphe 1.4, corollaire 2), au moins dans le cas qui est le n\^otre o\`u $\End(A)=\Z$, que si les repr\'esentations
$\rho'_{\ell,A_i}:\Gal(\bar{K}/K(\mu_{\ell^{\infty}}))\rightarrow \Aut(\T_{\ell}(A_i))$ sont isomorphes, alors $A_1$ et $A_2$ sont $\bar{K}$-isog\'enes. Le cas g\'en\'eral est sugg\'er\'e dans \cite{parzar} et trait\'e dans \cite{zar3}.

\medskip

\subsection{Point 1 du th\'eor\`eme \ref{prophodgeprod}}
 
Pour $i\in\{1,2\}$, notons $H_i=\hdg(A_i)=\Sp_{2g_i}$ et $H=\hdg(A_1\times A_2)$. On sait que $H\subset H_1\times H_2$ avec les projections $p_i:H\rightarrow H_i$ surjectives; on peut donc appliquer le lemme de Goursat (lemme \ref{goursat1}) : le groupe $N_1$ (resp. $N_2$) est un sous-groupe alg\'ebrique de $\Sp_{2g_1}$ (resp. $\Sp_{2g_2}$) et on dispose d'un isomorphisme $\rho:\Sp_{2g_1}/N_1\cong\Sp_{2g_2}/N_2$. Comme $\Sp_{2g_1}$ est presque simple on a $N_1=\{1\}$ , $\{\pm 1\}$ ou $\Sp_{2g_1}$. Dans le dernier cas, on voit que $N_2=\Sp_{2g_2}$ et donc $H=H_1\times H_2$ comme annonc\'e. Le cas $N_1=\{\pm 1\}$ est impossible car il imposerait 
$N_2=\{\pm 1\}$ (le groupe $\Sp_{2g}$ a pour centre $\{\pm 1\}$ alors que le centre de $\Sp_{2g}/\{\pm 1\}$ est trivial) et on obtiendrait un isomorphisme $\PSp_{2g_1}\rightarrow\PSp_{2g_2}$ qui impose $g_1=g_2$ comme tout automorphisme de $\PSp_{2g}$ provient d'un automorphisme  de $\Sp_{2g}$  (voir le lemme \ref{gsimple}) on en tirerait un $\alpha\in\Sp_{2g_1}(\Q)$ tel que
$H=\{(h_1,h_2)\in \Sp_{2g_1}\times\Sp_{2g_1}\;|\; h_2=\pm \alpha h_1\alpha^{-1}\}$ qui n'est pas connexe. Enfin le cas $N_1$ trivial entra\^{\i}ne $N_2$ trivial et l'isomorphisme 
$\rho:\Sp_{2g_1}\cong\Sp_{2g_2}$ impose $g_1=g_2$ et peut s'\'ecrire, comme tout automorphisme de $\Sp_{2g}$ est int\'erieur, $\rho(h)=\alpha h\alpha^{-1}$ avec $\alpha\in \Sp_{2g_1}(\Q)$. On obtiendrait donc

\begin{equation}
\hdg(A_1\times A_2)=\left\{  (h_1,h_2)\in \Sp_{2g_1}\times\Sp_{2g_1}\;|\; h_2=\alpha h_1 \alpha^{-1}\right\}
\end{equation}
mais alors un multiple entier de $\alpha$ induirait une isog\'enie entre $A_1$ et $A_2$.
$\findemo$

\subsection{Point 2 du th\'eor\`eme \ref{prophodgeprod}}

Nous allons prouver l'assertion \'equivalente concernant, non pas les groupes de Hodge, mais les groupes de Mumford-Tate : notons pour $i\in\{1,2\}$,
\[\rho_{\ell,i} : \Gal(\overline{K}/K)\rightarrow \Aut(\T_{\ell}(A_i))\subset \GSp_{2g_i}(\Z_{\ell})\text{ et }\psi_{\ell} : \Gal(\overline{K}/K)\rightarrow \Aut(\T_{\ell}(A_1))\times\Aut(\T_{\ell}(A_2))\]
\noindent les repr\'esentations $\ell$-adiques qui correspondent \`a l'action de Galois sur les points de torsion de $A_1$ et $A_2$ et o\`u $\psi_{\ell}=\rho_{\ell,1}\times\rho_{\ell,2}$. Rappelons que l'on note $\lambda_i : \GSp_{2g_i}\rightarrow \G_m$ l'application multiplicateur de noyau $\Sp_{2g_i}$. Notons enfin
\[E_{\ell}:=\left\{(x,y)\in \GSp_{2g_1}(\Z_{\ell})\times\GSp_{2g_2}(\Z_{\ell})\ |\ \lambda_1(x)=\lambda_2(y)\right\}\text{ et } H_{\ell}=\psi_{\ell}(\Gal(\overline{K}/K)).\]
Par hypoth\`eses, $\rho_{\ell,i}(\Gal(\overline{K}/K))= \GSp_{2g_i}(\Z_{\ell})$ \`a indice fini pr\`es, et il s'agit donc de montrer que $H_{\ell}$ est d'indice fini (d\'ependant \'eventuellement de $\ell$) dans $E_{\ell}$. Les groupes $H_{\ell}$ et $E_{\ell}$ sont des groupes de Lie $\ell$-adiques, il suffit donc de raisonner au niveau des alg\`ebres de Lie : notons $\hl$ et $\el$ les alg\`ebres de Lie respectives de $H_{\ell}$ et $E_{\ell}$ et prouvons que $\hl=\el$.

\begin{lemme}\label{alglie} Soient $i\in\{1,2\}$ et $g_1$, $g_2$ deux entiers strictement positifs. Notons $\lambda_i : \GSp_{2g_i}\rightarrow \G_m$ l'application multiplicateur de noyau $\Sp_{2g_i}$. Soit $\ell$ un nombre premier. Le groupe de Lie $\ell$-adique 
\[E_{\ell}:=\{(x,y)\in \GSp_{2g_1}(\Z_{\ell})\times\GSp_{2g_2}(\Z_{\ell})\ |\ \lambda_1(x)=\lambda_2(y)\}\]
\noindent a pour alg\`ebre de Lie l'alg\`ebre
\[\el=\{(x,y)\in \gsp_{2g_1}\times\gsp_{2g_2}\ |\ g_2\Tr(x)=g_1\Tr(y)\}.\]
\end{lemme}
\demo : Rappelons que, si $g\geq 1$ est un entier, on a l'identit\'e $\det=\lambda^g$ o\`u $\lambda$ est l'application multiplicateur : $\GSp_{2g} \rightarrow \G_m$. On sait que $\Lie(\det)=\Tr$ o\`u $\Tr$ d\'esigne l'op\'erateur trace. 
En passant aux applications tangentes, nous en d\'eduisons :
\[\Tr=g\Lie(\lambda).\]
\noindent Ceci nous permet d'obtenir l'alg\`ebre de Lie $\el$ sous la forme annonc\'ee.$\findemo$

\medskip

\noindent Par hypoth\`ese les projections $p_i : \hl\rightarrow \gsp_{2g_i}$ sont surjectives. On identifie $\gsp_{2g_1}$ avec $\gsp_{2g_1}\times\{0\}=\ker(p_2)$ et de m\^eme pour $\gsp_{2g_2}$ avec $\{0\}\times\gsp_{2g_2}=\ker(p_1)$. Posons $\n_i=\gsp_{2g_i}\cap \hl$. Par le lemme de Goursat (ou plut\^ot une variante \'evidente pour les alg\`ebres) l'image de $\hl$ dans $\gsp_{2g_1}/\n_1\times\gsp_{2g_2}/\n_2$ est le graphe d'un isomorphisme $\alpha : \gsp_{2g_1}/\n_1\rightarrow \gsp_{2g_2}/\n_2$. De plus par d\'efinition de $\hl$ et par le lemme \ref{alglie}, on sait que,
\[\n_1\subset\{(x,0)\ | \Tr(x)=0\}\subset \ssp_{2g_1}\text{ et } \n_2\subset\{(0,y)\ | \Tr(y)=0\}\subset \ssp_{2g_2}.\]
\noindent Or les alg\`ebres de Lie $\ssp_{2g_1}$ et $\ssp_{2g_2}$ sont simples, donc $\mathfrak{n}_{i}\in\{\{0\},\ssp_{2g_i}\}$.

\medskip

Si $\n_1=\ssp_{2g_1}$ alors, au vu de l'isomorphisme $\alpha$, on a n\'ecessairement $\n_2=\ssp_{2g_2}$ (et sym\'etriquement si $\n_2=\ssp_{2g_2}$ alors $\n_1=\ssp_{2g_1}$). Dans ce cas, on voit que 
\[\ssp_{2g_1}\times\ssp_{2g_2}\subset \hl\subset\el,\text{ et } \forall r\in\Q_{\ell}, \ r(I_{2g_1},I_{2g_2})\in\hl.\]
\noindent Soit donc $(x,y)\in\el$  et notons $r_0$ l'\'el\'ement de $\Q_{\ell}$ tel que $\Tr(x)=g_1r_0$. On a
\[\Tr(x)=\Tr(r_0 I_{2g_1})\text{ et }\Tr(y)=\frac{g_2}{g_1}\Tr(x)=g_2r_0=\Tr(r_0 I_{2g_2}).\]
\noindent Donc $(x,y)-(r_0 I_{2g_1},r_0 I_{2g_2})$ est dans $\hl$ donc $\el=\hl$.

\medskip

Sinon, $\n_1=\{0\}$ et $\n_2=\{0\}$ et $\alpha : \gsp_{2g_1}\rightarrow \gsp_{2g_2}$ est un isomorphisme. Montrons que ceci est impossible si $A_1$ n'est pas isog\`ene (sur $\overline{K}$) \`a $A_2$. Pour des raisons de dimension, l'isomorphisme $\alpha$ implique que $g_1=g_2$. Nous noterons d\'esormais $g$ cet entier.

\noindent On a
\[\gsp_{2g}=Z(\gsp_{2g})\oplus [\gsp_{2g},\gsp_{2g}],\text{ et } Z(\gsp_{2g})=\Q_{\ell}\cdot I_{2g}\text{ et }\ssp_{2g}=[\gsp_{2g},\gsp_{2g}].\]
\noindent De plus tout automorphisme (de $\gsp_{2g}$) respecte le centre et le groupe des commutateurs. Notamment, l'isomorphisme $\alpha$ envoie $\ssp_{2g}$ sur lui-m\^eme. Par ailleurs, on en d\'eduit \'egalement qu'il existe un \'el\'ement $r_0\in\Q_{\ell}$ tel que $\alpha(I_{2g})=r_0 I_{2g}$. Mais on sait par ailleurs que la compos\'ee $\lambda\circ \rho_{\ell} : \Gal(\overline{K}/K)\rightarrow \Z_{\ell}^{\times}$ n'est autre que le caract\`ere cyclotomique $\chi_{\text{cycl}}$. Vu la construction de $\alpha$, on en d\'eduit que $\Lie(\lambda)\circ \alpha=\Lie(\lambda)$. En \'evaluant cette identit\'e en $I_{2g}$, on conclut que $r_0=1$, autrement dit, $\alpha(I_{2g})=I_{2g}$ et donc la restriction de $\alpha$ \`a $\ssp_{2g}$ d\'etermine $\alpha$. Or tout automorphisme de $\mathfrak{sp}_{2g}$ est int\'erieur, i.e. il existe $f:V_{\ell}(A_1)\rightarrow V_{\ell}(A_2)$ qui est $\Q_{\ell}$-lin\'eaire telle que $\alpha_{|\ssp_{2g}}(u)=f\circ u\circ f^{-1}$ pour tout $u\in\ssp!
 _{2g}$. Soit $x\in\gsp_{2g}$. Notons $r\in\Q_{\ell}$ tel que $\Tr(x)=2rg$. On a $x-rI_{2g}\in\ssp_{2g}$ donc
\[\alpha(x)-rI_{2g}=\alpha(x-rI_{2g})=f\circ(x-rI_{2g})\circ f^{-1}=f\circ x\circ f^{-1}-rI_{2g}.\]
\noindent En particulier $\alpha(x)=f\circ x\circ f^{-1}$ pour tout $x\in\gsp_{2g}$, i.e. $\alpha$ est int\'erieur par $f$, donc on peut \'ecrire
$\hl$ comme l'ensemble des couples $(x,f\circ x\circ f^{-1})$ avec $x\in \gsp_{2g}$ et ainsi $f$ est un isomorphisme de $\hl$-modules. La th\'eorie de Lie montre alors qu'il existe un sous-groupe ouvert disons $U$ de $H_{\ell}$ tel que $f$ soit un isomorphisme de $U$-modules. Les repr\'esentations $\ell$-adiques $\rho_{\ell,1}$ et $\rho_{\ell,2}$ deviennent donc isomorphes sur l'extension $K'$ de $K$ correspondant \`a $U$;
ce qui , d'apr\`es les r\'esultats de Faltings (proposition \ref{lemfalt}), entra\^{\i}nerait que $A_1$ et $A_2$ sont isog\`enes et contredirait l'hypoth\`ese. $\findemo$

\subsection{Point 3 du th\'eor\`eme \ref{prophodgeprod}}

Reprenons les notations du paragraphe pr\'ec\'edent. Notons pour $i\in\{1,2\}$,
\[\rho_{\ell,i} : \Gal(\overline{K}/K)\rightarrow \Aut(\T_{\ell}(A_i))\subset \GSp_{2g_i}(\Z_{\ell})\text{ et }\psi_{\ell} : \Gal(\overline{K}/K)\rightarrow \Aut(\T_{\ell}(A_1))\times\Aut(\T_{\ell}(A_2))\]
\noindent les repr\'esentations $\ell$-adiques qui correspondent \`a l'action de Galois sur les points de torsion de $A_1$ et $A_2$ et o\`u $\psi_{\ell}=\rho_{\ell,1}\times\rho_{\ell,2}$. Rappelons que l'on note $\lambda_i : \GSp_{2g_i}\rightarrow \G_m$ l'application multiplicateur de noyau $\Sp_{2g_i}$. Notons enfin
\[E_{\ell}:=\left\{(x,y)\in \GSp_{2g_1}(\Z_{\ell})\times\GSp_{2g_2}(\Z_{\ell})\ |\ \lambda_1(x)=\lambda_2(y)\right\}\text{ et } H_{\ell}=\psi_{\ell}(\Gal(\overline{K}/K)).\]
Par hypoth\`eses, $\rho_{\ell,i}(\Gal(\overline{K}/K))= \GSp_{2g_i}(\Z_{\ell})$ \`a indice fini pr\`es ind\'ependant de $\ell$, et il s'agit de montrer que $H_{\ell}$ est d'indice fini born\'e ind\'ependamment de $\ell$ dans $E_{\ell}$.

\noindent Soit $H$ un sous-groupe ferm\'e d'indice born\'e par disons $C_0$ dans $\Sp_{2g}(\Z_{\ell})$. Sa projection est d'indice $\leq C_0$ dans $\Sp_{2g}(\F_{\ell})$ et donc, si $\ell$ est assez grand, est \'egale \`a $\Sp_{2g}(\F_{\ell})$ d'apr\`es les lemmes \ref{gsimple} et \ref{lemgs}. On conclut donc que $H=\Sp_{2g}(\Z_{\ell})$ d\`es que $\ell$ est assez grand  en utilisant le lemme~\ref{ssgsp}. 

\medskip

\noindent Introduisons les versions ``modulo $\ell$'' des objets pr\'ec\'edents : notons pour $i\in\{1,2\}$,
\[\bar{\rho}_{\ell,i} : \Gal(\overline{K}/K)\rightarrow \Aut(A_i[\ell])\subset \GSp_{2g_i}(\F_{\ell})\text{ et }\bar{\psi}_{\ell} : \Gal(\overline{K}/K)\rightarrow \Aut(A_1[\ell])\times\Aut(A_2[\ell])\]
\noindent les repr\'esentations $\ell$-adiques qui correspondent \`a l'action de Galois sur les points de $\ell$-torsion de $A_1$ et $A_2$ et o\`u $\bar{\psi}_{\ell}=\bar{\rho}_{\ell,1}\times\bar{\rho}_{\ell,2}$. Notons enfin
\[\bar{E}_{\ell}:=\left\{(x,y)\in \GSp_{2g_1}(\F_{\ell})\times\GSp_{2g_2}(\F_{\ell})\ |\ \lambda_1(x)=\lambda_2(y)\right\}\text{ et } \bar{H}_{\ell}=\bar{\psi}_{\ell}(\Gal(\overline{K}/K)).\]

\medskip

\noindent Nous allons tout d'abord prouver que si $\ell$ est assez grand, alors $\bar{H}_{\ell}=\bar{E}_{\ell}$. Nous conclurons ensuite.

\medskip

\noindent Notons $G_i:=\GSp_{2g_i}(\F_{\ell})$. Par hypoth\`ese les projections $p_i : \bar{H}_{\ell}\rightarrow G_i$ sont surjectives. On identifie $G_1$ avec $G_1\times\{1\}=\ker(p_2)$ et de m\^eme pour $G_2$ avec $\{1\}\times G_2=\ker(p_1)$. Posons $N_i=G_i\cap \bar{H}_{\ell}$. Par le lemme de Goursat l'image de $\bar{H}_{\ell}$ dans $G_1/N_1\times G_2/N_2$ est le graphe d'un isomorphisme $\alpha : G_1/N_1\rightarrow G_2/N_2$. De plus par d\'efinition de $\bar{H}_{\ell}$, on sait que, si $i\in\{1,2\}$,
\[N_1\subset\{(x,1)\ | \ \lambda_1(x)=1\}\subset \Sp_{2g_1}(\F_{\ell})\text{ et } N_2\subset\{(1,y)\ | \ \lambda_2(y)=1\}\subset \Sp_{2g_2}(\F_{\ell}).\]

\noindent Si $N_1=\Sp_{2g_1}(\F_{\ell})$, alors $N_2=\Sp_{2g_2}(\F_{\ell})$ (par cardinalit\'e). Dans ce cas, nous obtenons les inclusions
\[\Sp_{2g_1}(\F_{\ell})\times\Sp_{2g_2}(\F_{\ell})\subset \bar H_{\ell}\subset \bar E_{\ell}\subset \GSp_{2g_1}(\F_{\ell})\times\GSp_{2g_2}(\F_{\ell}).\]
\noindent On a de plus les projections surjectives 
\[\tilde{\lambda} : \GSp_{2g_1}(\F_{\ell})\times\GSp_{2g_2}(\F_{\ell})\overset{(\lambda_1,\lambda_2)}{\longrightarrow} \F_{\ell}^{\times}\times\F_{\ell}^{\times}\overset{\lambda_1\cdot\lambda_2^{-1}}{\longrightarrow}\F_{\ell}^{\times}\]

\begin{lemme}\label{groupeidiot}Soit $G_1,G_2,G_3$ trois groupes tels que $G_1\subset G_2$. Soit de plus $\phi : G_2\rightarrow G_3$ un morphisme de groupes tel que $\ker(\phi)\subset G_1$ et tel que $\phi(G_1)=\phi(G_2)$. Alors $G_1=G_2$.
\end{lemme}
\demo Imm\'ediat.$\findemo$

\medskip

\noindent Par d\'efinition $\bar{E}_{\ell}=\ker(\tilde{\lambda})$. Notons $\lambda_1'$ la restriction de $\lambda_1$ \`a $\bar{E}_{\ell}$. On voit alors que $\Sp_{2g_1}(\F_{\ell})\times\Sp_{2g_2}(\F_{\ell})=\ker(\lambda_1')$. De plus par d\'efinition de $\bar{H}_{\ell}$ et par la surjectivit\'e de $\bar{\rho}_{\ell,1}$ (si $\ell$ est assez grand), on sait que 
\[\lambda_1'(\bar{H}_{\ell})=\chi_{\textnormal{cycl}}(\Gal(\overline{K}/K))=\lambda_1'(\bar{E}_{\ell}),\]
\noindent o\`u $\chi_{\textnormal{cycl}}$ est le caract\`ere cyclotomique. En appliquant le lemme \ref{groupeidiot} avec $\phi=\lambda_1'$, $G_1=\bar{H}_{\ell}$, $G_2=\bar{E}_{\ell}$ et $G_3=\F_{\ell}^{\times}$, on conclut que $\bar{H}_{\ell}=\bar{E}_{\ell}$ comme d\'esir\'e.

\noindent Le groupe $N_1$ est distingu\'e dans $\GSp_{2g_1}(\F_{\ell})$, donc par le lemme \ref{gsimple}, si $N_1\not=\Sp_{2g_1}(\F_{\ell})$, alors $N_1\subset\{\pm 1\}$. Il en est alors de m\^eme pour $N_2$. Par cardinalit\'e on voit dans ce cas que $g_1=g_2$. Nous noterons $g$ cet entier dans la suite. De plus, le centre de $G_i/N_i$ est $\F_{\ell}^{\times}/N_i$ pour $i\in\{1,2\}$. L'isomorphisme $\alpha$ induit donc, par passage au quotient, un isomorphisme $\tilde{\alpha}$ de $G_1/\F_{\ell}^{\times}=\PGSp_{2g}(\F_{\ell})$ sur $\PGSp_{2g}(\F_{\ell})$.

\medskip

\noindent En appliquant le lemme \ref{pgsimple} on conclut qu'il existe un isomorphisme $f : A_1[\ell]\rightarrow A_2[\ell]$ tel que $\tilde{\alpha}(x)=f\circ x\circ f^{-1}$ pour tout $x\in\PGSp_{2g}(\F_{\ell})$. Ainsi, si $h=(x,y)\in\bar{H}_{\ell}$, il existe $\epsilon(h)\in\F_{\ell}^{\times}$ tel que 
\[y=\epsilon(h)f\circ x\circ f^{-1}.\]
\noindent En composant par l'application multiplicateur $\lambda$, on obtient $\epsilon(h)^2=1$ et on voit que $\epsilon$ est un homomorphisme de $\bar{H}_{\ell}$ dans $\{\pm1\}$. Posons $\epsilon_{\ell}=\epsilon\circ\bar{\psi}$. C'est un caract\`ere de $\Gal(\overline{K}/K)$ dans $\{\pm1\}$ tel que 
\[\forall x\in\Gal(\overline{K}/K),\ \ \bar{\rho}_{\ell,2}(x)=\epsilon_{\ell}(x)f\circ\bar{\rho}_{\ell,1}(x)\circ f^{-1}.\]

\medskip

\noindent Montrons maintenant que, si $\ell\geq 4g+1$, alors $\epsilon_{\ell}$ est non ramifi\'e en toute place ultram\'etrique non ramifi\'ee sur $\Q$ en laquelle $A_1$ et $A_2$ ont bonne r\'eduction. On suit pour cela l'argument de la preuve du lemme 8 de \cite{serreim72} en utilisant le corollaire 3.4.4. de \cite{raynaud} en lieu et place des corollaires 11 et 12 de \cite{serreim72} : soit $v$ une place ultram\'etrique du corps $K$ telle que $A_1$ et $A_2$ ont bonne r\'eduction en $v$ et que $v$ est non ramifi\'ee sur $\Q$. Supposons que la caract\'eristique de $v$ est $\ell$ (en effet $\epsilon_{\ell}$ est non-ramifi\'e en $v$ sinon \cite{serretate}) et notons $k_{\ell}$ une cl\^oture alg\'ebrique de $\F_{\ell}$. Notons par ailleurs $\chi_1,\ldots,\chi_{2g}$ (respectivement $\chi_1',\ldots,\chi_{2g}'$) les caract\`eres du groupe d'inertie mod\'er\'ee en $v$ \`a valeurs dans $k_{\ell}^{\times}$, intervenant dans le module galoisien $A_1[\ell]\otimes k_{\ell}$ (resp. $A_2[!
 \ell]\otimes k_{\ell}$), cf. \cite{serreim72} paragraphe 1.13. et \cite{raynaud} corollaire 3.4.4. En notant $\epsilon_v$ la restriction de $\epsilon_{\ell}$ au groupe d'inertie en $v$, on a pour tout $i$, (quitte \`a renum\'eroter les $\chi_i'$)
\[\chi_i=\epsilon_v\chi_i'\]
\noindent Comme l'indice de ramification $e(v)$ de $v$ est $1$, le corollaire 3.4.4. de \cite{raynaud} nous dit que les $\chi_i$ sont de la forme
\[\chi_i=\theta_{k_1}^{e(k_1)}\ldots\theta_{k_n}^{e(k_n)}\]
\noindent o\`u pour tout $r$, $e(r)\in\{0,1\}$ et o\`u les $\theta_{k_i}$ sont les $n$ caract\`eres fondamentaux de niveau $n$, l'entier $n$ pouvant varier dans $1,\ldots,2g$. Les invariants des $\chi_i$ et $\chi_i'$ dans $\Q/\Z$ (cf. \cite{serreim72} paragraphe 1.7) varient dans l'ensemble :
\[X=\left\{e(k_1)\frac{\ell^{k_1}}{\ell^n-1}+\cdots+e(k_n)\frac{\ell^{k_n}}{\ell^n-1}\ |\ k_i\in\{0,\ldots,n-1\},\ n\in\{1,\ldots,2g\}\right\}.\]
\noindent Enfin, comme $\epsilon_v^2=1$, son invariant est $0$ ou $\frac{1}{2}$ et est de la forme $x-x'$ avec $x,x'\in X$. Or si 
\[x=e(k_1)\frac{\ell^{k_1}}{\ell^n-1}+\cdots+e(k_n)\frac{\ell^{k_n}}{\ell^n-1},\]
\noindent 
\[0\leq x\leq \frac{n\ell^{n-1}}{\ell^n-1}<\frac{2g}{\ell-1}.\]
\noindent En particulier, $|x-x'|<\frac{2g}{\ell-1}$, et comme $\ell\geq 4g+1$, on voit que $|x-x'|<\frac{1}{2}$, donc n\'ecessairement l'invariant de $\epsilon_v$ vaut $0$, ce qui signifie que $\epsilon_v$ est non-ramifi\'e en $v$.

\medskip

\noindent Nous allons maintenant pouvoir conclure que $\bar{H}_{\ell}=\bar{E}_{\ell}$ pourvu que $\ell$ soit assez grand. Supposons par l'absurde qu'il existe une partie $L$ infinie de l'ensemble des nombres premiers telle que $\bar{H}_{\ell}\not=\bar{E}_{\ell}$ pour tout $\ell\in L$. Quitte \`a enlever un nombre fini de premiers, on peut supposer que $\ell\geq 4g+1$ et que les $\bar{\rho}_{\ell,i}$ sont surjectives sur $\Aut(A_i[\ell])$ pour tout $\ell\in L$. Soit $\ell\in L$. Les $\epsilon_{\ell}$ d\'efinis pr\'ec\'edemment sont non ramifi\'es en dehors d'un ensemble fini de places de $K$ ind\'ependant de $\ell$. Ceci implique que les $\epsilon_{\ell}$ varient dans un ensemble fini (quand $\ell$ est variable). Quitte \`a remplacer $L$ par une sous-partie infinie, on peut donc supposer que $\epsilon_{\ell}$ est ind\'ependant de $\ell$: notons le $\epsilon$. Sur l'extension $K'$ de $K$ correspondant au noyau de $\epsilon$, on obtient donc que les $\bar{\rho}_{\ell,i}$ sont i!
 somorphes. Par la seconde partie de la proposition \ref{lemfalt} ceci implique que $A_1$ et $A_2$ sont isog\`enes ce qui contredit l'hypoth\`ese.

\medskip

\noindent Conclusion : En invoquant le lemme \ref{relev},  on sait que $\bar{H}_{\ell}=\bar{E}_{\ell}$ si $\ell$ est assez grand. En particulier $\bar{H}_{\ell}$ contient $\Sp_{2g}(\F_{\ell})\times\Sp_{2g}(\F_{\ell})$ pour tout $\ell$ assez grand. Donc par le lemme \ref{relev} on conclut que $H_{\ell}$ contient $\Sp_{2g}(\Z_{\ell})\times\Sp_{2g}(\Z_{\ell})$ pour tout $\ell$ assez grand. De plus les projections de $H_{\ell}$ sur chacune des coordonn\'ees de $\GSp_{2g}(\Z_{\ell})\times\GSp_{2g}(\Z_{\ell})$ sont surjectives, donc $H_{\ell}=E_{\ell}$ pour tout $\ell$ assez grand, et en particulier, l'indice de $H_{\ell}$ dans $E_{\ell}$ est fini (par le point 2. pr\'ec\'edemment d\'emontr\'e) born\'e ind\'ependamment de $\ell$. $\findemo$

\section{Torsion pour un produit de vari\'et\'es ab\'eliennes de type $\GSp$}

Soit $d\geq 1$ un entier et pour tout $i$ compris entre $1$ et $d$, soient $n_i\geq 1$ des entiers et $A_i$ des vari\'et\'es ab\'eliennes de dimension respectives $g_i$, de type $\GSp$ et v\'erifiant la conjecture de Mumford-Tate. On note $A=\prod_{i=1}^dA_i^{n_i}$. Nous allons d\'emontrer dans ce paragraphe le th\'eor\`eme suivant.

\begin{theo}\label{th4}  Soit $A_i/K$ des vari\'et\'es ab\'eliennes non isog\`enes deux \`a deux, de dimension $g_i$, d\'efinie sur un corps de nombres, de type $\GSp$ et v\'erifiant la conjecture de Mumford-Tate. Soit 
$A:=A_1^{n_1}\times\dots\times A_r^{n_r}$ avec des entiers   $n_i\geq 1$. On a alors
\begin{equation}\label{valeurgammaproduitb}
\gamma(A)=\max_{\emptyset\not= I\subset\{1,\ldots,n\}}\left\{\frac{2\sum_{i\in I}n_i\dim A_i}{\dim\MT(\prod_{i\in I}A_i)}\right\}=\max_{\emptyset\not= I\subset\{1,\ldots,n\}}\left\{\frac{2\sum_{i\in I}n_ig_i}{1+\sum_{i\in I}2g_i^2+g_i}\right\}.
\end{equation}
\end{theo}

\subsection{Invariant $\gamma(A)$ pour un produit de vari\'et\'e ab\'eliennes de type \GSp}

Soient $\ell$ un nombre premier et $H$ un sous-groupe fini de $A[\ell^{\infty}]$. On peut facilement (voir par exemple \cite{hindry-ratazzi1} paragraphe 4.2) se ramener au cas o\`u $H$ est de la forme $H=\prod_{i=1}^dH_i^{n_i}$, les $H_i$ \'etant des sous-groupes finis de $A_i[\ell^{\infty}]$. Avec les notations introduites dans le cas simple (i.e. au paragraphe \ref{simple}), on peut, pour tout $i$, \'ecrire
\[H_i=\prod_{j=1}^{2g_i}\Z/\ell^{m_{ij}}\Z=\prod_{j=1}^{t_i}\left(\Z/\ell^{m(i)^j}\Z\right)^{\alpha_{ij}},\]
\noindent o\`u $(m(i)^j)_{j\geq 1}$ est une suite strictement d\'ecroissante. 
Par ailleurs, on introduit \'egalement pour tout $i$, une base $(\hat{e}_{ij})_j$ de $\T_{\ell}(A_i)$ associ\'ee \`a $H_i$ et on note comme pr\'ec\'edemment $G(H)$ et $G(H_i)$ les sous-groupes de $\Hdg(A)(\Z_{\ell})$(respectivement $\Hdg(A_i)(\Z_{\ell})$) stabilisant $H$ (respectivement $H_i$). Enfin on introduit \'egalement, comme dans le cas simple, pour tout $i$, la suite croissante de groupes alg\'ebriques sur $\Z_{\ell}$
\[\forall i\in\{1,\ldots,d\},\forall k\in\{1,\ldots,t_i\},\ \ G(H_i)_k:=\left\{M\in \Hdg(A_i)\ |\ \forall j\in I(i)_{t_i+1-k}\ \ M\hat{e}_{ij}=\hat{e}_{ij}\right\},\]
\noindent o\`u 
\[I(i)_r=\{j\in\{1,\ldots,2g_i\}\ |\ m_{ij}\geq m(i)^r\}\ \text{ est de cardinal }\ \sum_{j=1}^r\alpha_{ij}.\]
\noindent Enfin, on note, avec les notations du lemme \ref{calculH1}, $\delta:=\max\delta(H_i)$. Nous allons utiliser un r\'esultat galoisien que nous avons d\'emontr\'e dans \cite{hindry-ratazzi1} (voir th\'eor\`eme 6.6)~:

\begin{prop}\label{refauprodce}  Soient $A_1,\dots,A_r$ des vari\'et\'es ab\'eliennes d\'efinies sur $K$ et v\'erifiant la propri\'et\'e suivante~:  pour tout $i$ et tout groupe fini $H_i\subset A_i[\ell^{\infty}]$, il existe $m_i=m_i(H_i)$ tels qu'on a
\[K(H_i)\cap K(\mu_{\ell^{\infty}})\cong K(\mu_{\ell^{m_i}})\ ; \]
\noindent ainsi que l'identit\'e  o\`u $A:=A_1\times\cdots\times A_r$~:
\[\Gal(K(A[\ell^{\infty}])/K(\mu_{\ell^{\infty}}]))\cong \prod_{i=1}^r\Gal(K(A_i[\ell^{\infty}])/K(\mu_{\ell^{\infty}}]))\]
\noindent Alors si $m:=\max m_i $,   pour tout groupe fini $H=H_1\times\cdots\times H_r\subset A[\ell^{\infty}]$   on a $K(H)\cap K(\mu_{\ell^{\infty}})\cong K(\mu_{\ell^m})$ et
\[[K(H):K(\mu_{\ell^{m}})]\gg\ll\prod_{i=1}^{r}[K(H_i):K(\mu_{\ell^{m_i}})].\]
\end{prop}

\noindent Par le lemme \ref{calculH1}, on a
\[[K(H):K]=\delta(H)[\Hdg(A)(\Z_{\ell}):G(H)].\]
\noindent Or on sait par le paragraphe \ref{MTH} pr\'ec\'edent que dans notre situation $\Hdg(A)=\prod_{i=1}^d\Hdg(A_i)$ et que l'on a (\`a indice fini pr\`es)~:
\[\Gal(K(A[\ell^{\infty}])/K(\mu_{\ell^{\infty}}]))\cong \prod_{i=1}^r\Gal(K(A_i[\ell^{\infty}])/K(\mu_{\ell^{\infty}}]))\]
donc on peut appliquer la proposition \ref{refauprodce}. 
\[[K(H):K]=\delta(H)\prod_{i=1}^d[\Hdg(A_i)(\Z_{\ell}):G(H_i)].\]
\noindent Notons $d_{ik}$ la codimension  du groupe alg\'ebrique $G(H_i)_k$. Les calculs effectu\'es dans le cas simple nous donnent
\[\log_{\ell}[\Hdg(A)(\Z_{\ell}):G(H)]=\sum_{i=1}^d\sum_{k=1}^{t_i}d_{ik}\left(m(i)^{t_i+1-k}-m(i)^{t_i+1-(k-1)}\right).\]
\noindent De plus, quitte \`a renum\'eroter les $H_i$, on peut supposer (et on le fait) que $\delta(H)=\delta(H_1)$. On note alors $(\delta(1))_j$ la suite de $0$ et de $1$ relative \`a $\delta(H_1)$ d\'efinie dans le cas simple. On a
\[\log_{\ell}\delta\geq \sum_{j=1}^{t_1}\left(m(1)^{t_1+1-j}-m(1)^{t_1+1-(j-1)}\right)\delta(1)_j+O(1).\]
\noindent On pose par ailleurs $\delta(i)_j=0$ pour tout $j$ si $i\geq 2$. Avec ces notations, on trouve en suivant les calculs du cas simple,
\[\log_{\ell}[K(H):K]\geq\sum_{i=1}^d\sum_{j=1}^{t_i}m(i)^j\left[(\delta(i)_{t_i+1-j}-\delta(i)_{t_i+1-(j-1)})+(d_{i t_i+1-j}-d_{i t_i+1-(j-1)})\right]+O(1),\]
\noindent et
\[\log_{\ell}|H|=\sum_{i=1}^d\sum_{j=1}^{t_i}m(i)^jn_i\alpha_{ij}.\]
\noindent Notons 
\[b_{ij}:=(\delta(i)_{t_i+1-j}-\delta(i)_{t_i+1-(j-1)})+(d_{i t_i+1-j}-d_{i t_i+1-(j-1)}),\ \text{ et }\ a_{ij}:=n_i\alpha_{ij}.\]
\noindent On aura donc $|H|\ll [K(H):K]^{\gamma}$ si 
\[\gamma\geq\max\frac{\sum_{i=1}^d\sum_{j=1}^{t_i}m(i)^ja_{ij}}{\sum_{i=1}^d\sum_{j=1}^{t_i}b_{ij}},\]
\noindent le max \'etant pris sur les $m(i)^1\geq\ldots\geq m(i)^{t_i}$ pour $i$ entre $1$ et $d$. Comme dans le cas d'une vari\'et\'e ab\'elienne simple, on se ram\`ene alors au cas o\`u $H$ est d'exposant $\ell$. 

\medskip

\noindent Ainsi, en invoquant le lemme combinatoire \ref{combielem2} et en suivant les notations et calculs du cas simple, on voit que $|H|\ll [K(H):K]^{\gamma}$ si 
\[\gamma\geq \max \frac{\sum_{i=1}^d n_i(r(i)_{t_i+1-h_i}+s(i)_{t_i+1-h_i})}{\delta(1)_{t_1+1-h_1}+\sum_{i=1}^d\codim P_{r(i)_{t_i+1-h_i}, s(i)_{t_i+1-h_i}}}.\]    
\noindent Ce dernier max se r\'e\'ecrit sous la forme
\[\rho((n_i)_i,(g_i)_i):=\max_{0\leq s_i\leq r_i\leq g_i, 1\leq r_i} \frac{\sum_{i=1}^d n_i(r_i+s_i)}{\delta(1)+\sum_{i=1}^d (r_i+s_i)(2g_i-\frac{1}{2}(r_i+s_i-1)},\] \noindent et o\`u $\delta(1)=0$ si tout les $s_i$ sont nuls et $\delta(1)=1$ si l'un des $s_i$ est non nul. Il reste donc un calcul combinatoire pour conclure.

\subsection{Argument combinatoire}

\noindent On veut montrer que l'exposant
\[\alpha(A)=\max_{\emptyset\not= I\subset\{1,\ldots,d\}}\frac{2\sum_{i\in I}n_i\dim A_i}{\dim\MT\left(\prod_{i\in I}A_i\right)}\]
\noindent est admissible. Soit $n_i,g_i\geq 1$ pour $1\leq i\leq d$, on note $n=(n_1,\dots,n_d)$ et $g=(g_1,\dots,n_d)$. On d\'efinit la quantit\'e
\[\label{defalpha}
\alpha(n,g):=\max_{I\not=\emptyset}\left\{  \frac{2\sum_{i\in I}n_ig_i}{1+\sum_{i\in I}(2g_i^2+g_i)}\right\},
\]
o\`u le maximum est pris sur les sous-ensemble $I\subset\{1,\dots,d\}$. Dans notre situation $\alpha(A)=\alpha(n,g)$.

\noindent On pose \'egalement
\[\label{defrho}
\rho_0(n,g):=\max_{1\leq s_i\leq r_i\leq g_i}\left\{  \frac{\sum_{i=1}^dn_i(r_i+s_i)}{1+\sum_{i=1}^d(r_i+s_i)
(2g_i-\frac{1}{2}(r_i+s_i-1))}\right\},
\]
o\`u le maximum est pris sur les entiers $r_i$, $s_i$ tels que $1\leq s_i\leq r_i\leq g_i$. On peut clairement r\'ecrire ce dernier comme
\[\label{defrho2}
\rho_0(n,g):=\max_{2\leq x_i\leq 2g_i}\left\{  \frac{\sum_{i=1}^dn_ix_i}{1+\sum_{i=1}^dx_i
(2g_i-\frac{1}{2}(x_i-1))}\right\},
\]
Enfin on d\'efinit aussi
\[\label{defrho1}
\rho_1(n,g):=\max_{ r_i\leq g_i}\left\{  \frac{\sum_{i=1}^dn_ir_i}{\sum_{i=1}^dr_i
(2g_i-\frac{r_i-1}{2})}\right\},
\]
o\`u le maximum est pris sur les $r$-uplets  distincts de $(0,\dots,0)$. On a
\[\rho((n_i)_i,(g_i)_i)=\max(\rho_0(n,g),\rho_1(n_,g),\]
\noindent et on veut donc montrer le  r\'esultat suivant.

\begin{prop} On a l'in\'egalit\'e $\max\{\rho_1(n,g),\rho_0(n,g)\}\leq\alpha(A)$.
\end{prop}
\demo  On va utiliser l'in\'egalit\'e triviale suivante~:
\[\forall x\in[0,B],\; x^2-Bx\leq 0.\]
L'in\'egalit\'e
\[ \frac{\sum_{i=1}^dn_ix_i}{1+\sum_{i=1}^dx_i
(2g_i-\frac{1}{2}(x_i-1))}\leq\alpha\]
\'equivaut \`a 
\begin{equation}\label{ineqcombi}\sum_{i=1}^d\left(x_i^2-\left(4g_i+1-\frac{2n_i}{\alpha}\right)x_i\right)-2\leq 0.
\end{equation}
\noindent Observons que 
\[\alpha\geq\frac{2n_ig_i}{1+2g_i^2+g_i}=\frac{n_i}{g_i +\frac{1+g_i}{2g_i}}\geq \frac{n_i}{g_i+1}\] 
\noindent d'o\`u l'on tire 
\[2g_i-1\leq 4g_i+1-\frac{2n_i}{\alpha}.\]
\noindent Ainsi lorsque $x_i\leq 2g_i-1$ la contribution du terme en $x_i$ dans (\ref{ineqcombi}) est n\'egative et, si l'on note
$I:=\{i\in[1,d]\;|\; x_i=2g_i\}$ le membre de gauche de l'in\'egalit\'e (\ref{ineqcombi}) est major\'e par
\[\sum_{i\in I}\frac{4n_ig_i}{\alpha}-\sum_{i\in I}(4g_i^2+2g_i)-2\]
\noindent qui est n\'egatif si et seulement si
\[\alpha\geq \frac{\sum_{i\in I} 2n_ig_i}{1+\sum_{i\in I }(2g_i^2+g_i)}.\]
\noindent Le m\^eme calcul (en plus simple) montre que 
\[\rho_1(n,g)\leq \max_{I\not=\emptyset}\frac{\sum_{i\in I}g_in_i}{\sum_{i\in I}(g_i^2+g_i)}=
\max_{I\not=\emptyset}\frac{\sum_{i\in I}2g_in_i}{1+\sum_{i\in I}(2g_i^2+g_i)+\sum_{i\in I}g_i-1}\leq\alpha(n,g).\]
Ceci conclut.$\findemo$

\end{document}